\documentclass{article}
\usepackage{amsmath}
\usepackage{amssymb}
\usepackage{amsthm}
\usepackage[all]{xy}
\usepackage[latin1]{inputenc}        
\usepackage[dvips]{graphics}
\usepackage[dvips]{graphicx}
\usepackage{amscd}
\usepackage{mathrsfs}
\usepackage[mathcal]{eucal}
\usepackage{fullpage}
\usepackage{caption}
\usepackage{float}
\author{Ernesto Mistretta, Francesco Polizzi}

\title
{Standard isotrivial fibrations with $p_g=q=1$. II}
\date{}

\newtheorem{inizio}{Lemma}[section]
\newtheorem{theorem}[inizio]{Theorem}
\newtheorem{corollary}[inizio]{Corollary}
\newtheorem{proposition}[inizio]{Proposition}
\newtheorem{lemma}[inizio]{Lemma}

\newtheorem{definition}[inizio]{Definition}
\newtheorem*{teo}{Theorem}

\theoremstyle{definition}
\newtheorem{remark}[inizio]{Remark}

\newcommand{\g}{g_{\textrm{alb}}}
\newcommand{\lr}{\longrightarrow}

\newcommand{\mO}{\mathcal{O}}
\newcommand{\mZ}{\mathbb{Z}}
\newcommand{\sT}{\textrm{Sing}(T)}
\newcommand{\ST}{\emph{Sing}(T)}
\newcommand{\Sn}{\sum_{j=1}^s \left(1- \frac{1}{ \;n_j} \right)}

\newcommand{\sn}{1- \frac{1}{ \;n_1}}

\newcommand{\si}[2]{ \frac{1}{#1} (1,#2)}

\newcommand{\ssi}[6]{
\hline
$\si{#1}{#2}$ & $#3$ & $\si{#1}{#4}$  & $#5$ & $#6$ \\
}

\begin{document}


\maketitle

\abstract
A smooth, projective surface $S$ is called a
$\emph{standard isotrivial fibration}$ if there exists a finite group $G$ which acts
faithfully on two smooth projective curves $C$ and $F$ so that $S$ is isomorphic to
the minimal desingularization of $T:=(C \times F)/G$. Standard isotrivial fibrations
 of general type with
$p_g=q=1$ have been classified in \cite{Pol09} under the assumption
that $T$ has only Rational Double Points as singularities. In the
present paper we extend this result, classifying all cases where $S$
is a minimal model. As a by-product, we provide the first examples
of minimal surfaces of general type with $p_g=q=1$, $K_S^2=5$ and
Albanese fibration of genus $3$. Finally, we show with explicit
examples that the case where $S$ is not minimal actually occurs.

\endabstract


\section*{Introduction}
Surfaces of general type with $p_g=q=1$ are still not
well understood, and few examples are known. For a minimal surface
$S$ satisfying $p_g(S)=q(S)=1$, one has $2 \leq K_S^2 \leq 9$ and the Albanese map
is a connected fibration onto an elliptic curve. We denote by $\g$ the genus of a
general Albanese fibre of $S$. A classification of surfaces with $K_S^2=2,3$
 has been obtained by Catanese, Ciliberto, Pignatelli in
 \cite{Ca}, \cite{CaCi91}, \cite{CaCi93}, \cite{CaPi06}. For higher values of
$K_S^2$ some families are known, see \cite{Ca99}, \cite{CarPol},
\cite{Is05}, \cite{Pi08}, \cite{Pol09}, \cite{Ri07}, \cite{Ri08}. As
the title suggest, this paper considers surfaces with $p_g=q=1$
which are \emph{standard isotrivial fibrations}. This means that
there exists a finite group $G$ which acts faithfully on two smooth
projective curves $C$ and $F$ so that $S$ is isomorphic to the
minimal desingularization of $T:=(C \times F)/G$, where $G$ acts
diagonally on the product (see \cite{Se96}). When $p_g=q=1$ and $T$
contains at worst Rational Double Points (RDPs) as singularities,
standard isotrivial fibrations have been studied in \cite{Pol09} and
\cite{CarPol}. In the present article we make a further step toward
their complete classification, since we describe all cases where $S$
is a minimal model. As a by-product, we provide the first examples
of minimal surfaces of general type with
$p_g=q=1$, $K_S^2=5$ and $\g=3$ (see Section \ref{K2=5}). \\
Our classification procedure combines
methods from both geometry and combinatorial group theory. The basic idea is that since $S$
is the minimal desingularization of $T=(C \times F)/G$, the two projections $\pi_C \colon
C \times F \lr C$, $\pi_F \colon
C \times F \lr F$ induce two morphisms $\alpha \colon S \lr C/G$, $\beta \colon S \lr F/G$
 whose smooth fibres are isomorphic to $F$ and $C$, respectively. We have $1 =q(S)=g(C/G)
 +g(F/G)$, hence we may assume that $F/G \cong \mathbb{P}^1$ and $E:=C/G$ is an elliptic curve.
 Therefore $\alpha \colon S \lr E$ is the Albanese fibration of $S$ and $\g=g(F)$.
 The geometry of $S$ is encoded in the geometry of the two coverings $h \colon C \lr E$,
 $f \colon F \lr \mathbb{P}^1$ and the invariants of $S$ impose strong restrictions on $g(C)$,
 $g(F)$ and $|G|$. Indeed we can prove that under our assumptions $g(F)=2$ or $3$, hence we may
 exploit the classification of finite groups acting on curves of low genus given
 in $\cite{Br90}$. The problem of constructing our surfaces is then translated into the
 problem of finding two systems of generators of $G$, that we call $\mathcal{V}$ and
 $\mathcal{W}$, which are subject to strict conditions of combinatorial type.
The existence of such systems of generators can be checked in every case either by
hand-made computations or by using the computer algebra program \verb|GAP4| (see \cite{GAP4}). \\
This method of proof is similar to the one used in \cite{CarPol} and
\cite{Pol09}, of which the present paper is a natural sequel; the
main problem here is that when $T$ contains
 singularities worse than RDPs, they contribute not only to $\chi(\mO_S)$, but also to
$K_S^2$. However, since in any case $T$ contains only cyclic
 quotient singularities, this contribution is well known and can be computed in terms of
 Hirzebruch-Jung continued fractions (Corollary \ref{invariants-S}). When $S$ is minimal, we are able to use
 all this information in order to achieve a complete classification.
\begin{teo}
Let $\lambda \colon S \lr T:=(C \times F)/G$ be a standard isotrivial fibration
 of general type with $p_g=q=1$,
and assume that $T$ contains at least one
singularity which is not a \emph{RDP} and that $S$ is a minimal model.
Then there are exactly the following cases.
\begin{table}[ht!]
\begin{center}
\begin{tabular}{|c|c|c|c|c|c|}
\hline
$ $ & $g_{\textrm{alb}}=$ &  $ $ & $ $ &  \verb|IdSmall| & $ $ \\
$K_S^2$  & $g(F)$ & $g(C)$ & $G$ & \verb|Group|$(G)$ & Sing$(T)$  \\
\hline \hline
$5$ & $3$ & $3$ & $\mathcal{S}_3$ & $G(6,1)$ & $\frac{1}{3}(1,1)+\frac{1}{3}(1,2)$ \\
\hline
$5$ & $3$ & $5$ & $D_{4,3,-1}$ & $G(12,1)$ &  $\frac{1}{3}(1,1)+\frac{1}{3}(1,2)$ \\
\hline
$5$ & $3$ & $5$ & $D_6$ & $G(12,4)$ & $\frac{1}{3}(1,1)+\frac{1}{3}(1,2)$ \\
\hline
$5$ & $3$ & $9$ & $D_{2,12,5}$  & $G(24,5)$ & $\frac{1}{3}(1,1)+\frac{1}{3}(1,2)$ \\
\hline
$5$ & $3$ & $9$ & $\mathcal{S}_4$ & $G(24,12)$ & $\frac{1}{3}(1,1)+\frac{1}{3}(1,2)$ \\
\hline
$5$ & $3$ & $17$ & $\mZ_2 \times \mathcal{S}_4$ & $G(48,48)$ &
$\frac{1}{3}(1,1)+\frac{1}{3}(1,2)$ \\
\hline
$5$ & $3$ & $33$ & $\mathcal{S}_3 \ltimes (\mZ_4)^2$ & $G(96,64)$ &
$\frac{1}{3}(1,1)+\frac{1}{3}(1,2)$ \\
\hline
$5$ & $3$ & $57$ & $\textrm{PSL}_2(\mathbb{F}_7)$ &
$G(168,42)$ & $\frac{1}{3}(1,1)+\frac{1}{3}(1,2)$ \\
\hline
$3$ & $2$ & $11$ & $\mZ_2 \ltimes ((\mZ_2)^2 \times \mZ_3)$ &
$G(24,8)$ & $2 \times \frac{1}{2}(1,1)+\frac{1}{3}(1,1)+\frac{1}{3}(1,2)$\\
\hline
$3$ & $2$ & $21$ & $\textrm{GL}_2(\mathbb{F}_3)$ & $G(48,29)$ &
$2 \times \frac{1}{2}(1,1)+\frac{1}{3}(1,1)+\frac{1}{3}(1,2)$\\
\hline
$2$ & $2$ & $7$ & $D_{2,8,3}$ & $G(16,8)$ &
$2 \times \frac{1}{2}(1,1)+\frac{1}{4}(1,1)+\frac{1}{4}(1,3)$\\
\hline
$2$ & $2$ & $10$ & $\textrm{SL}_2(\mathbb{F}_3)$ & $G(24,3)$ &
$2 \times \frac{1}{2}(1,1)+\frac{1}{4}(1,1)+\frac{1}{4}(1,3)$\\
\hline
$2$ & $2$ & $3$ & $\mathcal{S}_3$ & $G(6,1)$ & $ 2 \times \frac{1}{3}(1,1)+
2 \times \frac{1}{3}(1,2)$ \\
\hline
$2$ & $2$ & $5$ & $D_{4,3,-1}$ & $G(12,1)$ & $ 2 \times \frac{1}{3}(1,1)+
2 \times \frac{1}{3}(1,2)$ \\
\hline
$2$ & $2$ & $5$ & $D_6$ & $G(12,4)$ & $ 2 \times \frac{1}{3}(1,1)+
2 \times \frac{1}{3}(1,2)$ \\
\hline
\end{tabular}
\end{center}
\end{table}
\end{teo}
Examples of non minimal standard isotrivial fibrations
  with $p_g=q=1$ actually exist. We exhibit two of them,
 one with $K_S^2=2$ (see Section \ref{sec-K2-2})
 and one with $K_S^2=1$ (see Section \ref{non-min-K2-1}); in both cases $\g=3$
 and the corresponding minimal model $\widehat{S}$ satisfies $K_{\widehat{S}}^2=3$.
The description of all non minimal examples would put an end to the classification
 of standard isotrivial fibrations with $p_g=q=1$; however, it seems to us difficult
 to achieve it by using our method. The main problem is that we
 are not able to find an effective lower bound for $K_S^2$. In fact, we can easily show that
 $S$ contains at most five $(-1)$-curves (Proposition \ref{-1-curves});
  nevertheless, when we contract them
  further $(-1)$-curves may appear. For instance, this actually happens
  in our example with
 $K_S^2=1$. \\ \\
$\mathbf{Notations \; and \; conventions}$. All varieties, morphisms, etc. in this
article are defined over $\mathbb{C}$. If $S$ is a projective,
non-singular surface $S$ then $K_S$ denotes the
canonical class, $p_g(S)=h^0(S, \; K_S)$ is the \emph{geometric genus}, $q(S)=h^1(S, \;
K_S)$ is the \emph{irregularity} and $\chi(\mathcal{O}_S)=1-q(S)+p_g(S)$ is the
\emph{Euler characteristic}. Throughout the paper we use the following notation for
groups:
\begin{itemize}
\item $\mZ_n$: cyclic group of order $n$.
\item $D_{p,q,r}=\mathbb{Z}_p \ltimes \mathbb{Z}_q= \langle x,y \; |
\; x^p=y^q=1, \; xyx^{-1}=y^r \rangle$: split metacyclic group of order $pq$. The group
$D_{2,n,-1}$ is the dihedral group of order $2n$ and it will be denoted by $D_n$.
\item $\mathcal{S}_n, \;\mathcal{A}_n$: symmetric, alternating group on $n$ symbols. We write the
composition of permutations from the right to the left; for instance, $(13)(12)=(123)$.
\item $\textrm{GL}_n(\mathbb{F}_q), \;\textrm{SL}_n(\mathbb{F}_q), \;
\textrm{PSL}_n(\mathbb{F}_ q)$: general linear, special linear and projective special linear
groups of $n
\times n$ matrices over a field with $q$ elements.
\item Whenever we give a presentation of a semi-direct product $H \ltimes N$, the
first generators represent $H$ and the last generators represent $N$. The action of $H$
on $N$ is specified by conjugation relations.
\item The order of a finite group $G$ is denoted by $|G|$. If $x \in G$, the order of
$x$ is denoted by $|x|$, its centralizer in $G$ by $C_G(x)$ and the conjugacy class of
$x$ by $\textrm{Cl}(x)$. If $x,y \in G$, their commutator is defined as
$[x,y]=xyx^{-1}y^{-1}$.
\item If $X= \{x_1,\ldots,x_n \} \subset G$, the subgroup generated
by $X$ is denoted by $\langle x_1,\ldots,x_n \rangle$. The derived subgroup of $G$ is
denoted by $[G,G]$.
\item \verb|IdSmallGroup|$(G)$ indicates the label of the group $G$ in
the  \verb|GAP4| database of small groups. For instance
\verb|IdSmallGroup|$(D_4)=G(8,3)$ means that $D_4$ is the third in the
 list of groups of order $8$.
\end{itemize}
\bigskip
$\mathbf{Acknowledgements.}$ The authors are indebted to I. Bauer,
F. Catanese and R. Pignatelli for many interesting discussions and
helpful suggestions. This research started at the University of
Bayreuth in September 2007, while the first author had a post-doc
position and the second author was visiting professor. Both authors
were supported by the \emph{DFG Forschergruppe ``Klassifikation
algebraischer Flächen und kompakter komplexer Mannigfaltigkeiten"}.

\section{Group-theoretic preliminaries}
In this section we fix the algebraic set-up and we present some preliminary results of
combinatorial type.
\begin{definition} \label{generating vect}
Let $G$ be a finite group and let
\begin{equation*}
\mathfrak{g}' \geq 0, \quad   m_r \geq m_{r-1} \geq \ldots \geq m_1 \geq 2
\end{equation*}
be integers. A \emph{generating vector} for $G$ of type $(\mathfrak{g}' \; | \; m_1,
\ldots ,m_r)$ is a $(2 \mathfrak{g}'+r)$-tuple of elements
\begin{equation*}
\mathcal{V}=\{g_1, \ldots, g_r; \; h_1, \ldots, h_{2\mathfrak{g}'} \}
\end{equation*}
such that the following conditions are satisfied:
\begin{itemize}
\item the set $\mathcal{V}$ generates $G$;
\item $|g_i|=m_i$;
\item $g_1g_2\cdots g_r \Pi_{i=1}^{\mathfrak{g}'} [h_i,h_{i+\mathfrak{g}'}]=1$.
\end{itemize}
If such a $\mathcal{V}$ exists, then $G$ is said to be $(\mathfrak{g}' \; | \; m_1,
\ldots ,m_r)-$\emph{generated}.
\end{definition}
\begin{remark} \label{abelian-gen}
If an \emph{abelian} group $G$ is $(\mathfrak{g}' \; | \; m_1,
\ldots ,m_r)-$generated then either $r=0$ or $r \geq 2$. Moreover if $r=2$ then
$m_1=m_2$.
\end{remark}
For convenience we make abbreviations such as $(4 \;| \; 2^3, 3^2)$ for $(4 \; | \;
2,2,2,3,3)$ when we write down the type of the generating vector $\mathcal{V}$.
 Moreover we set $\mathbf{m}:=(m_1, \ldots, m_r)$.
\begin{proposition}[Riemann Existence Theorem] \label{riemann ext}
A finite group $G$ acts as a group of automorphisms of some compact Riemann surface $X$
of genus $\mathfrak{g}$ if and only if there exist integers $\mathfrak{g}' \geq 0$ and
$m_r \geq m_{r-1} \geq \ldots \geq m_1 \geq 2$ such that $G$ is $(\mathfrak{g}'\; |\;
m_1, \ldots, m_r)-$generated, with generating vector $\mathcal{V}=\{g_1, \ldots, g_r; \;
h_1, \ldots, h_{2\mathfrak{g}'} \}$, and the Riemann-Hurwitz relation holds:
\begin{equation} \label{riemanhur}
2\mathfrak{g}-2=|G| \left( 2\mathfrak{g}'-2+\sum_{i=1}^r\bigg(1-\frac{1}{\;m_i} \bigg)
 \right).
\end{equation}
If this is the case, $\mathfrak{g}'$ is the genus of the quotient Riemann surface
$Y:=X/G$ and the $G-$cover $X \lr Y$ is branched in $r$ points $P_1, \ldots, P_r$ with
branching numbers $m_1, \ldots, m_r$, respectively. In addition, the subgroups $\langle
g_i \rangle$ and their conjugates provide all the nontrivial stabilizers of the action of
$G$ on $X$.
\end{proposition}
We refer the reader to [Br90, Section 2], [Bre00, Chapter 3],
\cite{H71} and [Pol09, Section 1] for more details. \\
Now let $X$ be a compact Riemann surface of genus $\mathfrak{g} \geq 2$ and
let $G \subseteq \textrm{Aut}(X)$.
For any $h \in G$ set $H:=\langle h \rangle$ and define the set of fixed points
of $h$ as
\begin{equation*}
\textrm{Fix}_X(h)=\textrm{Fix}_X(H):=\{x \in X \; |\; hx=x \}.
\end{equation*}
For our purposes it is also important to take into account how an
automorphism acts in a neighborhood of each of its fixed points. We
follow the exposition of \cite[pp.17, 38]{Bre00}. Let $\mathcal{D}$
be the unit disk and $h \in \textrm{Aut}(X)$ of order $m >1$ such that
$hx =x$ for a point $x \in X$.
 Then there is a unique
 primitive complex $m$-th root of unity $\xi$ such that any lift of
 $h$ to $\mathcal{D}$ that fixes a point in $\mathcal{D}$ is
 conjugate to the transformation $z \lr \xi \cdot z$ in
 $\textrm{Aut}(\mathcal{D})$. We write $\xi_x(h)=\xi$ and we
 call $\xi^{-1}$ the \emph{rotation constant} of $h$ in $x$.
Then for each integer $q \leq m-1$ such that $(q,m)=1$ we define
\begin{equation*}
\textrm{Fix}_{X,q}(h)=\{x \in \textrm{Fix}_X(h) \; | \;
\xi_x(h)=\xi^{q} \},
\end{equation*}
that is the set of fixed points of $h$ with rotation constant
$\xi^{-q}$. Clearly, we have
\begin{equation*}
\textrm{Fix}_X(h)=\biguplus_{\substack{q \leq m-1 \\ (q,m)=1}}\textrm{Fix}_{X,q}(h).
\end{equation*}
\begin{proposition} \label{fixed-points}
Assuming that we are in the situation of Proposition \emph{\ref{riemann ext}},
let $h \in G^{\times}$ be of order $m$, $H=\langle h \rangle$ and $(q,m)=1$.
Then
\begin{equation*} \label{formula-per-fix-tot}
|\emph{Fix}_{X}(h)|=|N_G(H)| \cdot
\sum_{\substack{1 \leq i \leq r \\ m|m_i
\\ H \; \sim_G \;  \langle g_i^{m_i/m} \rangle }} \frac{1}{\; m_i}
\end{equation*}

and

\begin{equation*} \label{formula-per-fix}
|\emph{Fix}_{X,q}(h)|=|C_G(h)| \cdot \sum_{\substack{1 \leq i \leq r \\ m|m_i
\\ h \; \sim_G \;  g_i^{m_iq/m} }} \frac{1}{\; m_i} ~.
\end{equation*}

\end{proposition}

\begin{proof}
See \cite[Lemma 10.4 and 11.5]{Bre00}.
\end{proof}
\begin{corollary} \label{conjugacy-q}
Assume that $h \sim_G h^q$. Then $|\emph{Fix}_{X,1}(h)|=|\emph{Fix}_{X,q}(h)|$.
\end{corollary}

\section{Surface cyclic quotient singularities and \\ Hirzebruch-Jung  resolutions}
\label{HJ-res}
Let $n$ and $q$ be natural numbers with $0 < q < n$, $(n,q)=1$
 and let $\xi_n$  be a primitive $n-$th root of unity.
Let us consider the action of the cyclic group $\mathbb{Z}_n=\langle \xi_n \rangle$ on
$\mathbb{C}^2$  defined by $\xi_n \cdot (x,y)=(\xi_nx, \xi_n^qy)$. Then the analytic
space $X_{n,q}=\mathbb{C}^2 / \mathbb{Z}_n$ has a cyclic quotient singularity of type
$\frac{1}{n}(1,q)$, and $X_{n,q} \cong X_{n', q'}$ if and only if $n=n'$ and either
$q=q'$ or $qq' \equiv 1$ (mod $n$). The exceptional divisor on the minimal
resolution $\tilde{X}_{n,q}$ of $X_{n,q}$ is a H-J string
 (abbreviation of Hirzebruch-Jung string), that is to say, a
 connected union $E=\bigcup_{i=1}^k Z_i$ of smooth rational curves $Z_1, \ldots, Z_k$ with
 self-intersection $\leq -2$, and ordered linearly so that $Z_i
 Z_{i+1}=1$ for all $i$, and $Z_iZ_j=0$ if $|i-j| \geq 2$.
More precisely, given the continued fraction
\begin{equation*}
\frac{n}{q}=[b_1,\ldots,b_k]=b_1-
                                \cfrac{1}{b_2 -\cfrac{1}{\dotsb
                                 - \cfrac{1}{\,b_k}}}, \quad b_i\geq 2 ~,
\end{equation*}
the dual graph of $E$ is  {\setlength{\unitlength}{1.1cm}
\begin{center}
\begin{picture}(1,0.5)
\put(0,0){\circle*{0.2}}
\put(1,0){\circle*{0.2}} \put(0,0){\line(1,0){1}}
\put(-0.3,0.2){\scriptsize $-b_1$}
\put(0.70,0.2){\scriptsize $-b_2$}
\put(2,0){\circle*{0.2}}
\put(1,0){\line(1,0){0.2}}
\put(1.3,0){\line(1,0){0.15}}
\put(1.55,0){\line(1,0){0.15}}
\put(1.8,0){\line(1,0){0.2}}
\put(3,0){\circle*{0.2}} \put(2,0){\line(1,0){1}}
\put(1.70,0.2){\scriptsize $-b_{k-1}$}
\put(2.70,0.2){\scriptsize $-b_k$}
\end{picture}         \hspace{2.5cm}
\end{center}
}

\vspace{.5cm}
\noindent (see \cite[Chapter II]{Lau71}). Notice that a RDP of type $A_n$ corresponds to the cyclic quotient singularity
$\frac{1}{n+1}(1,n)$.

\begin{definition} \label{numbers}
Let $x$ be a cyclic quotient singularity of type
$\frac{1}{n}(1,q)$. Then we set
\begin{equation*}
\begin{split}
h_x&=2- \frac{2+q+q'}{n}-\sum_{i=1}^k (b_i-2), \\
e_x&=k+1-\frac{1}{n}, \\
B_x&= 2e_x - h_x = \frac{1}{n} (q + q') + \sum_{i=1}^k b_i,
\end{split}
\end{equation*}
where $1\leq q' \leq n-1$ is such that $qq' \equiv 1$ $($\emph{mod}
$n)$.
\end{definition}

\section{Standard isotrivial fibrations}
In this section we establish the basic properties of standard isotrivial fibrations.
Definition \ref{def-stand} and Theorem \ref{Serrano}
 can be found in \cite{Se96}.
\begin{definition} \label{def-stand}
We say that a projective surface $S$ is
a  \emph{standard isotrivial fibration} if there exists a finite group
$G$ acting faithfully on two smooth projective curves $C$ and $F$ so that $S$ is
isomorphic to the minimal desingularization of $T:=(C \times F)/G$, where $G$ acts
diagonally on the product. The two maps $\alpha
\colon S \lr C/G$, $\beta \colon S \lr F/G$ will be referred as the \emph{natural
projections}.
\end{definition}
The stabilizer $H  \subseteq G$ of a point $y \in F$ is a cyclic group (\cite{FK92},
p.106). If $H$ acts freely on $C$, then $T$ is smooth along the scheme-theoretic fibre of
$\sigma \colon T \lr F/G$
 over $\bar{y} \in F/G$, and this fibre consists of the curve $C/H$
 counted with multiplicity $|H|$. Thus, the smooth fibres of $\sigma$
 are all isomorphic to $C$. On the contrary, if $ x \in C$ is fixed
 by some non-zero element of $H$, then $T$ has a cyclic quotient
 singularity  over the point $\overline{(x,y)} \in (C \times F)/G$.
These observations lead to the following
 statement, which describes the singular fibres that can arise in a
 standard isotrivial fibration (see \cite{Se96}, Theorem 2.1).
\begin{theorem} \label{Serrano}
Let $\lambda \colon S \lr T=(C \times F)/G$ be a standard isotrivial fibration and let us
consider the natural projection
 $\beta \colon S \lr F/G$.
Take any point over $\bar{y} \in F/G$ and let $\Lambda$ denote the schematic
fibre of $\beta$ over
$\bar{y}$. Then
\begin{itemize}
\item[$(i)$] The reduced structure of $\Lambda$ is the union of an
irreducible curve $Y$, called the central component of $\Lambda$,
 and either none or at least two mutually disjoint H-J strings, each
 meeting $Y$ at one point, and each being contracted by
 $\lambda$ to a singular point of $T$.
 These strings are in one-to-one
 correspondence with the branch points of $C \lr C/H$, where $H
 \subseteq G$ is the stabilizer of $y$.
\item[$(ii)$] The intersection of a string with $Y$ is transversal,
and it takes place at only one of the end components of the string.
\item[$(iii)$] $Y$ is isomorphic to $C/H$, and has multiplicity
equal to $|H|$ in $\Lambda$.
\end{itemize}
An analogous statement holds if we consider the natural projection
$\alpha \colon S \lr C/G$.
\end{theorem}

\begin{corollary} \label{signatures}
If $T$ has just two singularities, i.e.
\begin{equation*}
\emph{Sing}(T)=\frac{1}{n_1}(1,q_1)+\frac{1}{n_2}(1, q_2)
\end{equation*}
then $n_1=n_2$. \\ \\
If  $T$ has just three singularities, i.e.
\begin{equation*}
\emph{Sing}(T)=\frac{1}{n_1}(1,q_1)+\frac{1}{n_2}(1,q_2)+\frac{1}{n_3}(1,q_3)
\end{equation*}
then, for all $i=1, \,2, \,3$, the integer $n_i$ divides
$\emph{l.c.m.}\{n_k \, | \, k \neq i \}$
\end{corollary}
\begin{proposition} \label{all-sing-points}
Let $\lambda \colon S \lr T=(C \times F)/G$ be a standard isotrivial fibration. Assume
that
\begin{itemize}
\item[$(1)$] all elements of order $n$ are conjugate in $G$;
\item[$(2)$] $T$ contains a singular point of type $\frac{1}{n}(1,q)$ for some $q$ such
that $(q,n)=1$.
\end{itemize}
Then $T$ contains a singular point of type $\frac{1}{n}(1,r)$ for \emph{all} $r$ such
that $(r,n)=1$.
\end{proposition}
\begin{proof}
By assumption $(2)$ there exists a point $p=(p_1, p_2) \in C \times F$ such that the stabilizer of
$p$ has order $n$ and its generator $h$ acts, in suitable local coordinates centered at
$p$, as $h \cdot (x,\,y)=(\xi x, \, \xi^q y)$, where $\xi=e^{2 \pi i /n}$.  Therefore $p_2
\in |\textrm{Fix}_{F,q}(h)|$. Now let $r$ be such that $(r,n)=1$; using assumption
$(1)$ and Corollary \ref{conjugacy-q} we obtain
$|\textrm{Fix}_{F,r}(h)|=|\textrm{Fix}_{F,q}(h)| \neq 0$. If
$p_2' \in \textrm{Fix}_{F,r}(h)$, then in suitable local coordinates centered in
$p':=(p_1, \, p_2')$ the element $h$ acts as
 $h \cdot (x, \, y')=(\xi x, \, \xi^r y')$. This means that the image of $p'$
 in $T$ is a singular point of type $\frac{1}{n}(1,r)$.
\end{proof}
For a proof of the following result, see
\cite[p. 509-510]{Bar99} and \cite{Fre71}:
\begin{proposition}
Let $V$ be a smooth algebraic surface, and let $G$ be a finite group
acting on $V$ with only isolated fixed points. Let $\lambda \colon S
\lr T$ be the minimal desingularization. Then we have
\begin{itemize}
\item[$(i)$] $K_S^2 =\frac{1}{|G|} \cdot K_V^2 + \sum \limits_{x \in \emph{Sing}\; T}h_x$.
\item[$(ii)$] $e(S)=\frac{1}{|G|} \cdot e(V)+\sum \limits_{x \in \emph{Sing}\; T}e_x$.
\item[$(iii)$] $H^0(S, \Omega^1_S)=H^0(V, \Omega^1_V)^G$.
\end{itemize}
\end{proposition}
So we obtain
\begin{corollary} \label{invariants-S}
Let $\lambda \colon S \lr T=(C \times F)/G$ be a standard isotrivial
fibration. Then the invariants of $S$ are given by
\begin{itemize}
\item[$(i)$] $K_S^2 =\frac{8(g(C)-1)(g(F)-1)}{|G|} + \sum \limits_{x \in \emph{Sing}\; T}h_x$.
\item[$(ii)$] $e(S)=\frac{4(g(C)-1)(g(F)-1)}{|G|}+\sum \limits_{x \in \emph{Sing}\; T}e_x$.
\item[$(iii)$] $q(S)=g(C/G)+g(F/G)$.
\end{itemize}
\end{corollary}

\begin{remark}
If $g(C/G) >0$ and $g(F/G) >0$ then $S$ is necessarily a
minimal model. If instead $g(F/G)=0$ [respectively $g(C/G)=0$] it may happen
that the central component of some
reducible fibre of $\alpha$ [respectively $\beta$] is a
$(-1)$-curve. Examples of this situation
are given in Sections \ref{sec-K2-2} and \ref{non-min-K2-1}.
\end{remark}

\section{The case $\chi(\mO_S)=1$ }

\begin{proposition} \label{singularities}
Let  $\lambda \colon S \lr T=(C \times F)/G$ be a standard
isotrivial fibration with $\chi(\mO_S)=1$ and $K_S^2
\geq 2$. Then the possible singularities of $T$ are included in the
following list:

\vspace{.2cm}

\noindent
\begin{minipage}[t]{8cm}
$\bullet \;K_S^2=6:$
\begin{enumerate}
\item $2 \times \si{2}{1}.$
\end{enumerate}

\smallskip

$\bullet \;K_S^2=5:$
\begin{enumerate}
\item $ \si{3}{1} +  \si{3}{2};$
\item $2 \times \si{4}{1};$
\item $3 \times \si{2}{1}.$
\end{enumerate}

\smallskip

$\bullet \;K_S^2=4:$
\begin{enumerate}
\item $\si{4}{1} + \si{4}{3};$
\item $2 \times \si{5}{2};$
\item $\si{2}{1} + 2 \times \si{4}{1};$
\item $4 \times \si{2}{1}.$
\end{enumerate}

\smallskip

$\bullet \;K_S^2=3:$
\begin{enumerate}
\item $2 \times \si{4}{3};$
\item $\si{5}{1} + \si{5}{4};$
\item $\si{7}{2} + \si{7}{3};$
\item $\si{8}{1} + \si{8}{3};$
\item $\si{8}{5} + \si{8}{3};$
\item $\si{2}{1} + \si{4}{1}+ \si{4}{3};$
\item $2 \times \si{2}{1} + 2 \times \si{4}{1};$
\item $2 \times \si{2}{1} + \si{3}{1} + \si{3}{2};$
\item $5 \times \si{2}{1}.$
\end{enumerate}

\end{minipage}
\begin{minipage}[t]{8cm}

$\bullet \;K_S^2=2:$
\begin{enumerate}
\item $\si{6}{1} + \si{6}{5};$
\item $\si{9}{2} + \si {9}{4};$
\item $2 \times \si{10}{3};$
\item $\si{11}{3} + \si{11}{7};$
\item $\si{12}{5} + \si{12}{7};$
\item $2 \times \si{13}{5};$
\item $\si{2}{1} + 2 \times \si{4}{3};$
\item $\si{2}{1} + \si{5}{2} + \si{10}{3};$
\item $\si{2}{1} + \si{8}{1} + \si{8}{3};$
\item $\si{2}{1} + \si{8}{3} + \si{8}{5};$
\item $\si{3}{2} + 2 \times \si{6}{1};$
\item $\si{4}{1} + 2 \times \si{8}{3};$
\item $3 \times \si{5}{2};$
\item $2 \times \si{2}{1} +\si{4}{1} + \si{4}{3};$
\item $2 \times \si{2}{1} + 2 \times \si{5}{2};$
\item $4 \times \si{4}{1};$
\item $\si{3}{1}+\si{3}{2}+2 \times \si{4}{1};$
\item $2 \times \si{3}{1}+ 2 \times \si{3}{2};$
\item $3 \times \si{2}{1} + 2 \times \si{4}{1};$
\item $3 \times \si{2}{1} + \si{3}{1} +\si{3}{2};$
\item $6 \times \si{2}{1}.$
\end{enumerate}

\end{minipage}


 \vspace{.2cm}
Moreover the case $K_S^2 = 8$ occurs if and only if the action of $G$
is free, \emph{i.e.} if and only if $T$ is non-singular, whereas the
case $K_S^2 = 7$ does not occur.
\end{proposition}

\begin{proof}
By Corollary \ref{invariants-S} we have
$K_S^2= 2 e(S) - \sum_{x \in \textrm{Sing}\; T} (2e_x - h_x)$ and
Noether formula yields $K_S^2 = 12 - e(S)$, hence
\begin{equation}
\label{ksing}
K_S^2 = 8 - \frac{1}{3} \sum_{x \in \textrm{Sing}\; T}
B_x,
\end{equation}
where $B_x$ is as in Definition \ref{numbers}.

Notice that $3 \leq B_x \leq 18$ and that $B_x=3$ if and only if $x$
is of type $\frac{1}{2}(1,1)$. By Theorem \ref{Serrano} there are
either none or at least 2 singularities, and if there are exactly
two singularities they are of the form $\frac{1}{n}(1,q_1)$ and
$\frac{1}{n}(1, q_2)$, see Corollary \ref{signatures}. By analyzing
all  singularities with $B_x \leq 6$, we see that one cannot have
exactly two singularities $x_1$ and $x_2$ with $B_{x_1} > 12$ and
$B_{x_2} < 6$. Hence we may only consider singularities with $B_{x}
\leq 12$. A list of all such singularities with their numerical
invariants is given in Appendix $A$. \\For each fixed $K_S^2$ we
have to consider all possibilities for $\textrm{Sing}(T)$ such that
$\sum_{x \in \textrm{Sing}(T)}B_x=24-3K_S^2$ and we must exclude
those sets of singularities contradicting Corollary
\ref{signatures}. In this way we get our list. If $K_S^2=8$ then
equation (\ref{ksing}) implies that $T$ is smooth, whereas if
$K_S^2=7$ then $T$ would have exactly one singular point of type
$\frac{1}{2}(1,1)$, impossible by Theorem \ref{Serrano}.
\end{proof}

\begin{proposition} \label{divides}
Let $S$ be as in Proposition \emph{\ref{singularities}} and let us
assume $|\emph{Sing}\,T|=2$ or $3$. Then
\begin{itemize}
\item $m_i\;$ divides $\;g(C)-1\;$ for all $i \in \{1, \ldots,
r\}$, except at most one; \item $n_j\;$ divides $\;g(F)-1\;$ for
all $j \in \{1, \ldots, s \}$, except at most one.
\end{itemize}
If \,$|\emph{Sing}\,T|=4$ or $5$ then the same
statement holds with ``at most two'' instead of ``at most one''.
\end{proposition}
\begin{proof}
Assume $|\textrm{Sing}\,T|=2$ or $3$. Then by Theorem \ref{Serrano}
the corresponding H-J strings must belong to the same fibre of
$\beta \colon S \lr F/G$. It follows that, for all $i$ except one,
there is a subgroup $H$ of $G$, isomorphic to $\mathbb{Z}_{m_i}$,
which acts freely on $C$. Now Riemann-Hurwitz formula applied to
$C \lr C/H$ gives
\begin{equation*}
g(C)-1=m_i(g(C/H)-1),
\end{equation*}
 so $m_i$ divides $g(C)-1$. The statement about the $n_j$ is
 analogous. If $|\textrm{Sing}\,T|=4$ or $5$ then the H-J strings
 belong to at most two different fibres of $\beta$ and the same
 proof applies.
\end{proof}

\begin{corollary} \label{cor-divides}
If $\,|\emph{Sing}\,T| \leq 3$ and $g(F)=2$ then $s=1$, that is
$\mathbf{n}=(n_1)$. In particular, under these assumptions $G$ is
not abelian $($see \emph{Remark} \emph{\ref{abelian-gen}}$)$.
\end{corollary}

\section{Standard isotrivial fibrations with $p_g=q=1$}

From now on we suppose that $\lambda \colon S \lr T=(C \times F)/G$
is a standard isotrivial fibration with $p_g=q=1$. Since $q=1$, we
may assume that $E:=C/G$ is an elliptic curve and that $F/G \cong
\mathbb{P}^1$.
Then the natural projection $\alpha \colon S \lr E$ is the Albanese
morphism of $S$ and $g_{\textrm{alb}}=g(F)$. Let $\mathcal{V}=\{g_1,
\ldots g_r \}$ be a generating vector for $G$ of type $(0\; | \;
m_1, \ldots , m_r)$, inducing the $G$-cover $F \lr \mathbb{P}^1$ and
let $\mathcal{W}=\{\ell_1, \ldots \ell_s; \; h_1, h_2\}$ be a
generating vector of type $(1\; | \; n_1, \ldots , n_s)$ inducing the $G$-cover
$C\lr E$. Then Riemann-Hurwitz formula implies
\begin{equation} \label{generi}
\begin{split}
2g(F)-2 &=|G|\bigg(-2+\sum_{i=1}^r \bigg( 1- \frac{1}{\;m_i} \bigg)
\bigg) \\
2g(C)-2 & =|G| \sum_{j=1}^s \bigg(1-\frac{1}{\;n_j} \bigg).
\end{split}
\end{equation}
The cyclic subgroups $\langle g_1 \rangle,
\ldots ,\langle g_r \rangle$ and their conjugates provide the
non-trivial stabilizers of the action of $G$ on $F$, whereas
$\langle \ell_1 \rangle, \ldots, \langle \ell_s \rangle$ and their
conjugates provide the non-trivial stabilizers of the actions of $G$
on $C$. The singularities of $T$ arise from the points in $C \times
F$ with nontrivial stabilizer; since the action of $G$ on $C \times
F$ is the diagonal one, it follows that the set
 $\mathscr{S}$ of all nontrivial stabilizers for the action of $G$ on $C \times F$
 is given by
\begin{equation} \label{stabilizzatori}
\mathscr{S}= \bigg( \bigcup_{\sigma \in G} \bigcup_{i=1}^r \langle
\sigma g_i \sigma^{-1} \rangle \bigg) \cap \bigg( \bigcup_{\sigma
\in G} \bigcup_{j=1}^s \langle \sigma \ell_j \sigma^{-1} \rangle
\bigg) \cap G^{\times}.
\end{equation}

\begin{proposition} \label{strukture}
Let $G$ be a finite group which is both $(0 \; |\;m_1, \ldots,
m_r)-$generated and \\$(1\; | \;n_1, \ldots, n_s)-$generated, with
generating vectors $\mathcal{V}=\{g_1, \ldots, g_r \}$ and
$\mathcal{W}=\{\ell_1, \ldots, \ell_s; \; h_1,h_2 \}$, respectively.
Denote by
\begin{equation*}
\begin{split}
f & \colon F \lr \mathbb{P}^1=F/G, \\   h & \colon C \lr E= C/G
\end{split}
\end{equation*}
the $G$-covers induced by $\mathcal{V}$ and $\mathcal{W}$ and
let $g(F), \; g(C)$ be the genera of $F$ and $C$, that are related
to $|G|, \; \mathbf{m}, \; \mathbf{n}$ by \emph{(\ref{generi})}.
Define
\begin{equation*}
k=\frac{8(g(C)-1)(g(F)-1)}{|G|}+\sum_{x \in \emph{Sing}(T)}h_x
\end{equation*}
and assume that equality
\begin{equation} \label{ksing-1}
k=8 - \frac{1}{3} \sum_{x \in \emph{Sing}(T)}B_x
\end{equation}
holds. Then the minimal desingularization $S$ of $T$ satisfies
\begin{equation*}
p_g(S)=q(S)=1, \quad  K_S^2=k.
\end{equation*}
Moreover, if $k>0$ then $S$ is of general type.
\end{proposition}
\begin{proof}
The normal surface $T$ satisfies $q(T)=1$; since all quotient singularities are rational
it follows $q(S)=1$. Corollary \ref{invariants-S} and
relation \eqref{ksing-1} yield $K_S^2=k$ and $K_S^2+e(S)=12$, hence $\chi(\mO_S)=1$
by Noether formula; this implies $p_g(S)=1$. Finally if $k >0$ then $S$
 is of general type, because $q(S) >0$.
\end{proof}

\begin{lemma} \label{fund-rel}
Let $\lambda \colon S \lr T=(C \times F)/G$ be a standard isotrivial
fibration with $p_g=q=1$. Then we have
\begin{equation} \label{eq:fund-form}
K_S^2-\sum_{x \in \emph{Sing}(T)}h_x=4(g(F)-1) \Sn.
\end{equation}
\end{lemma}
\begin{proof}
Applying Corollary \ref{invariants-S} and the second
relation in \eqref{generi} we obtain
\begin{equation*}
\begin{split}
K_S^2-\sum_{x \in \textrm{Sing}(T)}h_x&=4(g(F)-1) \cdot 2
\frac{(g(C)-1)}{|G|} \\ &=4(g(F)-1) \Sn.
\end{split}
\end{equation*}
\end{proof}

The cases where $T$ has at worst RDP have already been classified in
\cite{Pol09} and \cite{CarPol}. Hence, in the sequel we will
consider the situation where $T$ contains at least one singularity
which is not a RDP.

\begin{proposition} \label{possibilities}
Let $\lambda \colon S \lr T=(C \times F)/G$ be a standard
isotrivial fibration with $p_g=q=1$, $K_S^2 \geq 2$ such that $T$ contains at least
one singularity which is not a \emph{RDP}. Then there are at most the
following possibilities: \\ \\
$\bullet$ $K_S^2=5$
\begin{itemize}
\item[] $g(F)=3, \quad \mathbf{n}=(3), \;\;\quad \ST=\frac{1}{3}(1,1)+
  \frac{1}{3}(1,2);$
\item[] $g(F)=3, \quad \mathbf{n}=(8), \;\;\quad \ST= 2\times \frac{1}{4}(1,1).$
\item[]
\end{itemize}
$\bullet$ $K_S^2=4$
\begin{itemize}
\item[] $g(F)=3, \; \;\quad \mathbf{n}=(4), \quad \ST=\frac{1}{2}(1,1)+2 \times \frac{1}{4}(1,1)$.
\item[]
\end{itemize}
$\bullet$ $K_S^2=3$
\begin{itemize}
\item[] $g(F)=2, \quad \mathbf{n}=(2,4), \quad \ST=2 \times\frac{1}{2}(1,1)+2 \times \frac{1}{4}(1,1);$
\item[] $g(F)=2, \quad \mathbf{n}=(6), \quad \quad \ST=2 \times\frac{1}{2}(1,1)+ \frac{1}{3}(1,1)+\frac{1}{3}(1,2)$.
\item[]
\end{itemize}
$\bullet$ $K_S^2=2$
\begin{itemize}
\item[] $g(F)=3, \quad \mathbf{n}=(16), \quad \; \ST=\frac{1}{2}(1,1)+ \frac{1}{8}(1,1)+\frac{1}{8}(1,3);$
\item[] $g(F)=2, \quad \mathbf{n}=(8), \; \;\quad \;\ST=\frac{1}{2}(1,1)+ \frac{1}{8}(1,3)+\frac{1}{8}(1,5);$
\item[] $g(F)=3, \quad \mathbf{n}=(12), \quad \; \ST=\frac{1}{3}(1,2)+ 2 \times \frac{1}{6}(1,1);$
\item[] $g(F)=2, \quad \mathbf{n}=(4),  \;  \; \quad \; \ST= 2\times\frac{1}{2}(1,1)+ \frac{1}{4}(1,1)+\frac{1}{4}(1,3);$
\item[] $g(F)=3, \quad \mathbf{n}=(4), \; \; \;  \quad  \ST= 4 \times \frac{1}{4}(1,1);$
\item[] $g(F)=2, \quad \mathbf{n}=(4^2), \,\,\quad \ST =4 \times \frac{1}{4}(1,1);$
\item[] $g(F)=2, \quad \mathbf{n}=(3), \, \, \, \, \, \quad \ST= 2 \times \frac{1}{3}(1,1)
+ 2 \times \frac{1}{3}(1,2)$.
\end{itemize}
\end{proposition}
\begin{proof}
For every value of $K_S^2$ we must analyze all possible
singularities of $T$ as listed in Proposition \ref{singularities}.
Moreover we have to exclude the cases in which all singularities of
$T$ are RDPs, namely $K_S^2=6$, $K_S^2 = 5 ~ (iii)$, $K_S^2 = 4 ~
(iv)$,
 $K_S^2 = 3 ~(ix)$ and $K_S^2 = 2 ~ (xxi)$, where $T$ contains only
 singular points of type $A_1$, and $K_S^2 = 3 ~ (i)$, where $T$ contains
 only singular points of type $A_3$.
 \\
\\
$\bullet$ $K_S^2=5$ \\ \\
$(i). \; \; \sT=\frac{1}{3}(1,1)+\frac{1}{3}(1,2)$. Using formula
(\ref{eq:fund-form}) and the table in Appendix A  we obtain
\begin{equation*}
(g(F)-1) \Sn =\frac{4}{3}.
\end{equation*}
If $s=1$ then $\frac{4}{3} < g(F)-1 \leq \frac{8}{3}$, which implies
$g(F)=3$,  $\mathbf{n}=(3)$. If $s \geq 2$ then $g(F)-1 \leq \frac{4}{3}$, so $g(F)=2$ which
contradicts Corollary \ref{cor-divides}. \\ \\
$(ii). \; \;  \sT=2 \times \frac{1}{4}(1,1)$. We obtain
\begin{equation*}
(g(F)-1) \Sn =\frac{7}{4}.
\end{equation*}
If $s=1$ then $\frac{7}{4} < g(F)-1 \leq \frac{7}{2}$, hence $g(F)=3$ or $4$. The case
$g(F)=4$ is numerically impossible, so the only possibility is
$g(F)=3$, $\mathbf{n}=(8)$.
If $s \geq2$ then $g(F)-1 \leq \frac{7}{4}$, so $g(F)=2$ which
contradicts Corollary \ref{cor-divides}. \\ \\
$\bullet$ $K_S^2=4$ \\ \\
$(i). \; \; \sT=\frac{1}{4}(1,1) +  \frac{1}{4}(1,3)$. We obtain
\begin{equation*}
(g(F)-1) \Sn =\frac{5}{4}.
\end{equation*}
If $s=1$ then $\frac{5}{4} < g(F)-1 \leq \frac{5}{2}$, so $g(F)=3$ which is
impossible. If $s \geq 2$ then $g(F)-1 \leq \frac{5}{4}$, so $g(F)=2$
 which contradicts Corollary \ref{cor-divides}. \\ \\
$(ii). \; \; \sT=2 \times \frac{1}{5}(1,2)$. We obtain
\begin{equation*}
(g(F)-1) \Sn =\frac{6}{5}.
\end{equation*}
If $s=1$ then $\frac{6}{5} < g(F)-1 \leq \frac{12}{5}$, so $g(F)=3$ which is
impossible. If $s \geq 2$ then $g(F)-1 \leq \frac{6}{5}$, so
$g(F)=2$ which contradicts Corollary \ref{cor-divides}. \\ \\
$(iii). \; \;  \sT=\frac{1}{2}(1,1)+2 \times \frac{1}{4}(1,1)$. We
obtain
\begin{equation*}
(g(F)-1) \Sn =\frac{3}{2}.
\end{equation*}
If $s=1$ then $\frac{3}{2} < g(F)-1 \leq 3$, so
$g(F)=3$ or $4$. In the former case we obtain the possibility
$g(F)=3$, $\mathbf{n}=(4)$; in the latter $\mathbf{n}=(2)$ and $T$ would
contain only singular points of of type $A_1$, a contradiction.
If $s \geq 2$ then $g(F)-1 \leq \frac{3}{2}$, so $g(F)=2$ against
 Corollary \ref{cor-divides} \\ \\
$\bullet$ $K_S^2=3$ \\ \\
$(ii). \; \; \sT= \frac{1}{5}(1,1) + \frac{1}{5}(1,4)$. We obtain
\begin{equation*}
(g(F)-1)\Sn=\frac{6}{5}.
\end{equation*}
If $s=1$ then $\frac{6}{5} < g(F)-1 \leq \frac{12}{5}$, so $g(F)=3$
which is numerically impossible. If $s \geq 2$ then $g(F)-1 \leq \frac{6}{5}$,
 so $g(F)=2$ that contradicts
Corollary \ref{cor-divides}. \\ \\
$(iii). \; \; \sT= \frac{1}{7}(1,2)+ \frac{1}{7}(1,3)$. We obtain
\begin{equation*}
(g(F)-1)\Sn=\frac{8}{7},
\end{equation*}
which gives $g(F)-1 \leq \frac{16}{7}$, so either $g(F)=2$ or
$g(F)=3$. In the former case we must have $s=1$, which is
impossible. In the latter we obtain $\Sn=\frac{4}{7}$, which
has no integer solutions. \\ \\
$(iv). \; \; \sT= \frac{1}{8}(1,1)+\frac{1}{8}(1,3)$. We obtain
\begin{equation*}
(g(F)-1)\Sn=\frac{17}{8},
\end{equation*}
which implies $g(F)-1 \leq \frac{17}{4}$, hence $2 \leq g(F) \leq 5$.
If $g(F)=2$ then $s=1$ by Corollary \ref{cor-divides}, and this is
numerically impossible.
It follows $g(F) =3, 4$ or $5$, hence
$\Sn=\frac{17}{16}, \, \frac{17}{24}$ or $\frac{17}{32}$, respectively. In
all cases there are no solutions. \\ \\
$(v). \; \; \sT=\frac{1}{8}(1,5)+\frac{1}{8}(1,3)$. We obtain
\begin{equation*}
(g(F)-1)\Sn=\frac{9}{8},
\end{equation*}
which implies either $g(F)=2$ or $g(F)=3$. In the former case we have
$s =1$, which is  numerically impossible. In the latter we obtain
$\Sn=\frac{9}{16}$, which has no solutions. \\ \\
$(vi). \; \; \sT=\frac{1}{2}(1,1)+ \frac{1}{4}(1,1)+
\frac{1}{4}(1,3)$. We obtain
\begin{equation*}
(g(F)-1)\Sn=1,
\end{equation*}
hence either $g(F)=2$ or $g(F)=3$. The former case yields
$\mathbf{n}=(2^2)$ and the latter $\mathbf{n}=(2)$; then $T$ would
have at worst $A_1-$singularities, a contradiction. \\ \\
$(vii). \; \; \sT= 2 \times \frac{1}{2}(1,1) + 2 \times
\frac{1}{4}(1,1)$. We obtain
\begin{equation*}
(g(F)-1)\Sn=\frac{5}{4},
\end{equation*}
hence either $g(F)=2$ or $g(F)=3$. In the former case the only possibility
is $\mathbf{n}=(2,4)$. In the latter we obtain $\Sn=\frac{5}{8}$,
 which has no solutions. \\ \\
$(viii). \; \;  \sT=2 \times \frac{1}{2}(1,1) +
\frac{1}{3}(1,1)+\frac{1}{3}(1,2)$. We obtain
\begin{equation*}
(g(F)-1)\Sn=\frac{5}{6},
\end{equation*}
which gives the only possibility
$g(F)=2$, $\mathbf{n}=(6)$. \\ \\
$\bullet K^2=2$. \\ \\
$(i). \; \; \sT= \frac{1}{6}(1,1)+\frac{1}{6}(1,5)$. We obtain
\begin{equation*}
(g(F)-1)\Sn=\frac{7}{6},
\end{equation*}
hence $g(F)=3$ or $g(F)=2$. If $g(F)=3$ then $\Sn=\frac{7}{12}$, which
is a contradiction. If $g(F)=2$ then $s=1$ by Corollary \ref{cor-divides}, so
$\sn=\frac{7}{6}$ which is impossible.    \\ \\
$(ii). \; \; \sT= \frac{1}{9}(1,2)+ \frac{1}{9}(1,4)$. We have
\begin{equation*}
(g(F)-1)\Sn=\frac{10}{9},
\end{equation*}
hence $g(F)=2$ or $3$ and we obtain a contradiction as before.\\ \\
$(iii). \; \; \sT= 2 \times \frac{1}{10}(1,3)$. We have
\begin{equation*}
(g(F)-1)\Sn=\frac{11}{10},
\end{equation*}
hence $g(F)=2$ or $3$ and we obtain a contradiction as before. \\ \\
$(iv). \; \; \sT=\frac{1}{11}(1,3)+ \frac{1}{11}(1,7)$. We have
\begin{equation*}
(g(F)-1) \Sn=\frac{12}{11},
\end{equation*}
hence $g(F)=2$ or $3$ and we obtain a contradiction as before. \\ \\
$(v). \; \; \sT =\frac{1}{12}(1,5)+\frac{1}{12}(1,7)$. We have
\begin{equation*}
(g(F)-1)\Sn=\frac{13}{12},
\end{equation*}
hence $g(F)=2$ or $3$ and we obtain a contradiction as before. \\ \\
$(vi). \; \; \sT= 2 \times \frac{1}{13}(1,5)$. We have
\begin{equation*}
(g(F)-1)\Sn=\frac{14}{13},
\end{equation*}
hence $g(F)=2$ or $3$ and we obtain a contradiction as before. \\ \\
$(viii). \; \; \sT= \frac{1}{2}(1,1) +
\frac{1}{5}(1,2)+\frac{1}{10}(1,3)$. We have
\begin{equation*}
(g(F)-1)\Sn=\frac{9}{10},
\end{equation*}
hence $g(F)=2$ and $\mathbf{n}=(10)$. This means that $G$ is
one of the groups listed
in Table $2$ of Appendix $B$ and that $10$ divides $|G|$, a contradiction. \\ \\
$(ix). \; \; \sT= \frac{1}{2}(1,1)+
\frac{1}{8}(1,1)+\frac{1}{8}(1,3)$. We have
\begin{equation*}
(g(F)-1)\Sn=\frac{15}{8},
\end{equation*}
so $g(F)=2, \,3$ or $4$. If $g(F)=2$ then $s=1$ (Corollary \ref{cor-divides}),
which gives a contradiction. The case $g(F)=4$ is numerically impossible.
Finally, if $g(F)=3$ we obtain $\mathbf{n}=(16)$. \\ \\
$(x). \; \; \sT= \frac{1}{2}(1,1)+
\frac{1}{8}(1,3)+\frac{1}{8}(1,5)$. We obtain
\begin{equation*}
(g(F)-1)\Sn= \frac{7}{8}
\end{equation*}
and the only possibility is $g(F)=2$, $\mathbf{n}=(8)$. \\ \\
$(xi). \; \; \sT= \frac{1}{3}(1,2) + 2 \times \frac{1}{6}(1,1)$. We
have
\begin{equation*}
(g(F)-1)\Sn= \frac{11}{6},
\end{equation*}
hence $g(F)=2, \,3$ or $4$. If $g(F)=2$ then $s=1$, a
contradiction. The case $g(F)=4$ is numerically impossible. Finally, if
$g(F)=3$ we obtain $\mathbf{n}=(12)$.  \\ \\
$(xii). \; \; \sT= \frac{1}{4}(1,1)+ 2 \times \frac{1}{8}(1,3)$. We
obtain
\begin{equation*}
(g(F)-1)\Sn = \frac{5}{4},
\end{equation*}
hence $g(F)=3$ or $2$. If $g(F)=3$ then $\Sn=\frac{5}{8}$, which has no solutions; if
$g(F)=2$ then $s=1$ which is
a contradiction. \\ \\
$(xiii). \; \; \sT= 3 \times \frac{1}{5}(1,2)$. We obtain
\begin{equation*}
(g(F)-1)\Sn= \frac{4}{5},
\end{equation*}
hence $\mathbf{n}=(5)$ and $g(F)=2$. This means that $5$ divides $|G|$ and that $G$
 is one of the groups listed
in Table $2$ of Appendix $B$, a contradiction. \\ \\
$(xiv) \; \; \sT= 2 \times \frac{1}{2}(1,1)+\frac{1}{4}(1,1)+
\frac{1}{4}(1,3)$. We obtain
\begin{equation*}
(g(F)-1)\Sn= \frac{3}{4},
\end{equation*}
which gives the possibility $g(F)=2$, $\mathbf{n}=(4)$. \\ \\
$(xv). \; \; \sT= 2 \times \frac{1}{2}(1,1)+ 2 \times
\frac{1}{5}(1,2)$. We obtain
\begin{equation*}
(g(F)-1)\Sn= \frac{7}{10},
\end{equation*}
which is impossible. \\ \\
$(xvi). \; \; \sT= 4 \times \frac{1}{4}(1,1)$. We obtain
\begin{equation*}
(g(F)-1)\Sn=\frac{3}{2},
\end{equation*}
so $g(F)=2, \, 3$ or $4$. At least one of the $n_i$ must be divisible by
$4$, otherwise $T$ could not contain singularities of type $\frac{1}{4}(1,1)$.
Hence the only possibilities are $g(F)=2$, $\mathbf{n}=(4^2)$
and  $g(F)=3$, $\mathbf{n}=(4)$. \\ \\
$(xvii). \; \; \sT=\si{3}{1}+\si{3}{2}+2 \times \si{4}{1}$. We obtain
\begin{equation*}
(g(F)-1)\Sn=\frac{13}{12},
\end{equation*}
so either $g(F)=2$ or $g(F)=3$. Consequently, either
$\Sn =\frac{13}{12}$ or $\Sn=\frac{13}{24}$,
and in both cases there are no solutions. \\ \\
$(xviii). \; \; \sT= 2 \times \si{3}{1} + 2 \times \si{3}{2}$. We obtain
\begin{equation*}
(g(F)-1)\Sn= \frac{2}{3},
\end{equation*}
which gives the possibility $g(F)=2$, $\mathbf{n}=(3)$. \\ \\
$(xix). \; \; \sT= 3 \times \frac{1}{2}(1,1) + 2 \times
\frac{1}{4}(1,1)$. We obtain
\begin{equation*}
(g(F)-1)\Sn=1,
\end{equation*}
hence $g(F)=2$ or $3$. If $g(F)=2$ then $\mathbf{n}=(2^2)$, whereas if $g(F)=3$
then $\mathbf{n}=(2)$; both cases are impossible otherwise $T$ would have only
$A_1$-singularities. \\ \\
$(xx). \; \; \sT= 3 \times \si{2}{1}+ \si{3}{1}+ \si{3}{2}$. We obtain
\begin{equation*}
(g(F)-1)\Sn=\frac{7}{12},
\end{equation*}
which has no solutions. \\ \\
This concludes the proof of Proposition \ref{possibilities}.
\end{proof}

\subsection{The case where $G$ is abelian}

\begin{proposition}
Let $\lambda \colon S \lr T=(C \times F)/G$ be a standard
isotrivial fibration with $p_g=q=1$, $K_S^2 \geq 2$ and $G$
abelian. Then $T$ contains at worst \emph{RDPs}.
\end{proposition}
\begin{proof}
Suppose that $G$ is abelian and $T$ contains at least one singularity
which is not a RDP.
Then by Proposition \ref{possibilities} and
Remark \ref{abelian-gen} we must have
\begin{equation*}
K_S^2=2, \quad g(F)=2, \quad \mathbf{n}=(4^2), \quad
 \textrm{Sing}(T) =4 \times \frac{1}{4}(1,1).
\end{equation*}
Corollary \ref{invariants-S} implies $g(C)-1=\frac{3}{4}|G|$.
Referring to Table 1 of Appendix $B$, we are left with two cases: \\ \\
$\bullet$ $(1c)$. $\;G=\mZ_4, \quad \mathbf{m}=(2^2, 4^2), \quad
g(C)=4$.\\
$\bullet$ $(1h)$. $\;G=\mZ_8, \quad \mathbf{m}=(2,8^2), \quad
\;  g(C)=7$.\\ \\
Therefore $G$ must be cyclic. Let $\mathcal{W}=\{\ell_1, \ell_2; \, h_1, h_2 \}$ be a
 generating vector of type $(1\, | \, 4^2)$ for $G$; then $\ell_1=(\ell_2)^{-1}$ and
 Proposition \ref{fixed-points} implies
\begin{equation*}
|\textrm{Fix}_{C,1}(\ell_1)|=|\textrm{Fix}_{C,3}(\ell_1)|=2.
\end{equation*}
In particular $\textrm{Fix}_{C,1}(\ell_1)$ and $\textrm{Fix}_{C,3}(\ell_1)$ are both
 nonempty. Hence the same argument used in proof of Proposition
\ref{all-sing-points} shows that if $S$ contains a singularity of type
 $\frac{1}{4}(1,1)$ then it contains also a
 singularity of type $\frac{1}{4}(1,3)$, a contradiction.
This concludes the proof.
\end{proof}
Therefore in the sequel we may assume that $G$ is a nonabelian group.

\subsection{The case $K_S^2 = 5$} \label{K2=5}

\begin{lemma} \label{no-8-gen}
Referring to Table $3$ in Appendix B, in cases
 $(3g)$, $(3h)$, $(3i)$, $(3j)$, $(3s)$, $(3u)$, $(3v)$
the group $G$ is not $(1\,|\,8)-$generated.
\end{lemma}
\begin{proof}
In cases $(3g)$, $(3h)$, $(3i)$, $(3j)$ we have
$[G,G]=\mZ_2$, in case $(3s)$ we have $[G,G]= \mathcal{A}_4$, in case $(3u)$
we have $[G,G]=\mZ_4 \times \mZ_4$ and in case $(3v)$ we have $[G,G]=G(48,3)$.
So in all cases $[G,G]$ contains no elements
of order $8$ and we are done.
\end{proof}

\begin{proposition}
Let $\lambda \colon S \lr T=(C \times F)/G$ be a standard
isotrivial fibration with $p_g=q=1$, $K_S^2=5$ such that
$T$ contains at least one singularity which is not a \emph{RDP}.
Then
\begin{equation*}
g(F)=3, \quad \mathbf{n}=(3), \quad \ST=\si{3}{1} + \si{3}{2}.
\end{equation*}
Furthermore exactly the following cases occur:
\begin{center}
\begin{tabular}{|c|c|c|c|c|}
\hline
$ $ & \verb|IdSmall| & $ $ & & \\
$G$ & \verb|Group|$(G)$ & $\mathbf{m}$ & $g(C)$ & \emph{Is} $S$ \emph{minimal?}\\
\hline \hline $\mathcal{S}_3$ & $G(6,1)$ & $(2^4,3)$ & $3$  & \emph{Yes}\\
\hline $D_{4,3,-1}$ & $G(12,1)$ & $(4^2 , 6)$ & $5$ & $\emph{Yes}$ \\
\hline $D_{6}$ & $G(12,4)$ & $(2^3 , 6)$ & $5$ & $\emph{Yes}$ \\
\hline $D_{2,12,5}$ & $G(24,5)$ & $(2 , 4, 12)$ & $9$ & $\emph{Yes}$
\\ \hline $\mathcal{S}_4$ & $G(24,12)$ & $(3 , 4^2)$ & $9$ & $\emph{Yes}$ \\ \hline
$\mathcal{S}_4$ & $G(24,12)$ & $(2^3 , 3)$ & $9$ & $\emph{Yes}$ \\
\hline $\mathbb{Z}_2 \times \mathcal{S}_4$  & $G(48,48)$ & $(2,4,6)$
& $17$ & $\emph{Yes}$
\\ \hline $\mathcal{S}_3 \ltimes (\mathbb{Z}_4)^2$ & $G(96,64)$ &
$(2, 3, 8)$ & $33$ & $\emph{Yes}$ \\ \hline
$\emph{PSL}_2(\mathbb{F}_7)$ & $G(168,42)$ & $(2,3,7)$  & $57$ & $\emph{Yes}$\\
\hline
\end{tabular}
\end{center}
\end{proposition}
\begin{proof}
If $K_S^2=5$, by Proposition \ref{possibilities} we have two possibilities:
\begin{equation*}
\begin{split}
& (a) \; \; g(F)=3, \quad \mathbf{n}=(3), \quad \sT= \frac{1}{3}(1,1)+ \frac{1}{3}(1,2); \\
& (b) \; \; g(F)=3, \quad \mathbf{n}=(8), \quad \sT=2 \times
\frac{1}{4}(1,1).
\end{split}
\end{equation*}
In particular $G$ must be one of the groups in Table $3$ of Appendix $B$.

Let us rule out first case $(b)$.
If it occurs then $(g(C)-1)=\frac{7}{16} |G|$ by Corollary \ref{invariants-S},
 so $|G|$ is divisible by $16$;
moreover, since $\mathbf{n}=(8)$, the group $G$ must be $(1\,|\,8)$-generated.
Cases $(3g)$, $(3h)$, $(3i)$, $(3j)$,
$(3s)$, $(3u)$, $(3v)$ are excluded by Lemma \ref{no-8-gen};
 cases $(3q)$, $(3r)$, $(3t)$ are excluded by
Proposition \ref{divides}.
So $(b)$ does not occur and we must only consider possibility $(a)$.

If it occurs then $g(C)-1 = \frac{1}{3}|G|$, so $|G|$ is divisible by $3$.
Moreover, since $\mathbf{n}=(3)$,
the group $G$ must be $(1\,|\,3)$-generated.
Cases $(3f)$, $(3k)$, $(3m)$, $(3n)$, $(3t)$, $(3u)$
 in Table $3$ are excluded by
Proposition \ref{divides}.
Now let us show that all the remaining cases occur. \\ \\

$\bullet$ Case $(3a)$. $\;G=\mathcal{S}_3, \; \;
\mathbf{m}=(2^4,3),\; \; g(C)=3,
\; \; \textrm{Sing}(T)=\frac{1}{3}(1,1)+\frac{1}{3}(1,2)$. Set
\begin{equation*}
\begin{split}
g_1 &=(12), \; \; g_2=(12), \; \; g_3=(12), \; \; g_4=(13), \; \; g_5=(123) \\
\ell_1 &=(123), \; \; h_1=(13), \; \; h_2=(12).
\end{split}
\end{equation*}
We have $\mathscr{S}= \textrm{Cl((123))}=
\{(123), \, (132)\}$ and for all $h \in \mathscr{S}$
\begin{equation*}
\begin{split}
|\textrm{Fix}_{F,1}(h)|&=|\textrm{Fix}_{F,2}(h)|=1\\
|\textrm{Fix}_{C,1}(h)|&=|\textrm{Fix}_{C,2}(h)|=1.
\end{split}
\end{equation*}
So $C \times F$ contains exactly four points with nontrivial stabilizer
 and for each of them the $G$-orbit has
cardinality $|G|/|(123)|=2$. Hence $T$ contains
 precisely two singular points and looking at the rotation constants we
see that $\textrm{Sing}(T)=
\frac{1}{3}(1,1)+\frac{1}{3}(1,2)$, as required. So this case occurs by
 Proposition \ref{strukture}. \\ \\

$\bullet$ Case $(3d)$. $\; G=D_{4,3,-1}= \langle x,y\; | \; x^4=y^3=1,
 \; xyx^{-1}=y^{-1} \rangle, \; \; \mathbf{m}=(4^2,6),\;
\; g(C)=5,
\; \; \textrm{Sing}(T)=\frac{1}{3}(1,1)+\frac{1}{3}(1,2)$. Set
\begin{equation*}
\begin{split}
g_1 &=x, \; \; g_2=xy, \; \; g_3=y^2x^2  \\
\ell_1 &=y, \; \; h_1=y, \; \; h_2=x.
\end{split}
\end{equation*}
We have $\mathscr{S}= \textrm{Cl}(y)=\{ y, y^2\}$ and for all $h \in \mathscr{S}$
\begin{equation*}
\begin{split}
|\textrm{Fix}_{F,1}(h)|&=|\textrm{Fix}_{F,2}(h)|=1 \\
|\textrm{Fix}_{C,1}(h)|&=|\textrm{Fix}_{C,2}(h)|=2.
\end{split}
\end{equation*}
So $C \times F$ contains exactly $8$ points with nontrivial
stabilizer and for each of them the  $G$-orbit has
cardinality $|G|/|y|=4$. Looking at the rotation constants we see
that $\textrm{Sing}(T)=
\frac{1}{3}(1,1)+\frac{1}{3}(1,2)$, as required. \\ \\

$\bullet$ Case $(3e)$. $\; G=D_6=\langle  x,y\; | \; x^2=y^6=1,
 \; xyx^{-1}=y^{-1} \rangle, \; \; \mathbf{m}=(2^3,6),\; \;
g(C)=5,
\; \; \textrm{Sing}(T)=\frac{1}{3}(1,1)+\frac{1}{3}(1,2)$. Set
\begin{equation*}
\begin{split}
g_1 &=x, \; \; g_2=xy^2, \; \; g_3=y^3, \; \; g_4=y  \\
\ell_1 &=y^2, \; \; h_1=x, \; \; h_2=y.
\end{split}
\end{equation*}
We have $\mathscr{S}= \textrm{Cl}(y^2)=\{ y^2, y^4\}$ and
for all $h \in \mathscr{S}$
\begin{equation*}
\begin{split}
|\textrm{Fix}_{F,1}(h)|&=|\textrm{Fix}_{F,2}(h)|=1 \\
|\textrm{Fix}_{C,1}(h)|&=|\textrm{Fix}_{C,2}(h)|=2.
\end{split}
\end{equation*}
So $C \times F$ contains exactly $8$ points  with nontrivial
stabilizer and for each of them the $G$-orbit has
cardinality $|G|/|y^2|=4$. Looking at the rotation constants we see
that $\textrm{Sing}(T)=
\frac{1}{3}(1,1)+\frac{1}{3}(1,2)$, as required. \\ \\

$\bullet$ Case $(3l)$. $\; G=D_{2,12,5}=\langle x,y\; | \; x^3=y^7=1,
 \; xyx^{-1}=y^2 \rangle, \; \;
\mathbf{m}=(2,4,12),\; \; g(C)=9,
\; \; \textrm{Sing}(T)=\frac{1}{3}(1,1)+\frac{1}{3}(1,2)$. Set
\begin{equation*}
\begin{split}
g_1 &=x, \; \; g_2=xy^{11}, \; \; g_3=y  \\
\ell_1 &=y^4, \; \; h_1=y, \; \; h_2=x.
\end{split}
\end{equation*}
We have $\mathscr{S}= \textrm{Cl}(y^4)= \{ y^4, y^8\}$ and
for all $h \in \mathscr{S}$
\begin{equation*}
\begin{split}
|\textrm{Fix}_{F,1}(h)|&=|\textrm{Fix}_{F,2}(h)|=1 \\
|\textrm{Fix}_{C,1}(h)|&=|\textrm{Fix}_{C,2}(h)|=4.
\end{split}
\end{equation*}
So $C \times F$ contains exactly $16$ points with nontrivial stabilizer and
for each of them the $G$-orbit has
cardinality $|G|/|y^4|=8$. Looking at the rotation constants we see
that $\textrm{Sing}(T)=
\frac{1}{3}(1,1)+\frac{1}{3}(1,2)$, as required. \\ \\

$\bullet$ Case $(3o)$. $\; G=\mathcal{S}_4, \; \;
\mathbf{m}=(3,4^2),\; \; g(C)=9,
\; \; \textrm{Sing}(T)=\frac{1}{3}(1,1)+\frac{1}{3}(1,2)$. Set
\begin{equation*}
\begin{split}
g_1 &=(123), \; \; g_2=(1234), \; \; g_3=(1243)  \\
\ell_1 &=(123), \; \; h_1=(142), \; \; h_2=(23).
\end{split}
\end{equation*}
We have $\mathscr{S}= \textrm{Cl}((123))$, hence $|\mathscr{S}|=8$
 and for all $h \in \mathscr{S}$
\begin{equation*}
\begin{split}
|\textrm{Fix}_{F,1}(h)|&=|\textrm{Fix}_{F,2}(h)|=1\\
|\textrm{Fix}_{C,1}(h)|&=|\textrm{Fix}_{C,2}(h)|=1.
\end{split}
\end{equation*}
So $C \times F$ contains exactly $16$ points with nontrivial stabilizer
and for each of them the $G$-orbit has
cardinality $|G|/|(123)|=8$. Looking at the rotation constants we
see that $\textrm{Sing}(T)=
\frac{1}{3}(1,1)+\frac{1}{3}(1,2)$, as required. \\ \\

$\bullet$ Case $(3p)$. $\; G=\mathcal{S}_4, \; \;
\mathbf{m}=(2^3,3),\; \; g(C)=9,
\; \; \textrm{Sing}(T)=\frac{1}{3}(1,1)+\frac{1}{3}(1,2)$. Set
\begin{equation*}
\begin{split}
g_1 &=(12), \; \; g_2=(24), \; \; g_3=(13)(24), \; \; g_4=(123)  \\
\ell_1 &=(123), \; \; h_1=(142), \; \; h_2=(23).
\end{split}
\end{equation*}
We have $\mathscr{S}= \textrm{Cl}((123))$, hence $|\mathscr{S}|=8$ and
for all $h \in \mathscr{S}$
\begin{equation*}
\begin{split}
|\textrm{Fix}_{F,1}(h)|&=|\textrm{Fix}_{F,2}(h)|=1\\
|\textrm{Fix}_{C,1}(h)|&=|\textrm{Fix}_{C,2}(h)|=1.
\end{split}
\end{equation*}
So $C \times F$ contains exactly $16$ points with nontrivial stabilizer
and for each of them the $G$-orbit has
cardinality $|G|/|(123)|=8$. Looking at the rotation constants we
see that $\textrm{Sing}(T)=
\frac{1}{3}(1,1)+\frac{1}{3}(1,2)$, as required. \\ \\

$\bullet$ Case $(3s)$. $\; G= \mathbb{Z}_2 \times \mathcal{S}_4, \;
\; \mathbf{m}=(2,4,6),\; \; g(C)=17,
\; \; \textrm{Sing}(T)=\frac{1}{3}(1,1)+\frac{1}{3}(1,2)$. \\
Put $\mathbb{Z}_2= \langle z\, | \, z^2=1 \rangle$ and set
\begin{equation*}
\begin{split}
g_1 &=z(14), \; \; g_2=(1234), \; \; g_3=z(132)  \\
\ell_1 &=(123), \; \; h_1=z(142), \; \; h_2=z(23).
\end{split}
\end{equation*}
We have $\mathscr{S}= \textrm{Cl}((123))$, hence $|\mathscr{S}|=8$
 and for all $h \in \mathscr{S}$
\begin{equation*}
\begin{split}
|\textrm{Fix}_{F,1}(h)|&=|\textrm{Fix}_{F,2}(h)|=1 \\
|\textrm{Fix}_{C,1}(h)|&=|\textrm{Fix}_{C,2}(h)|=2.
\end{split}
\end{equation*}
So $C \times F$ contains exactly $32$ points with nontrivial stabilizer and
for each of them the $G$-orbit has
cardinality $|G|/|(123)|=16$. Looking at the rotation constants we
see that $\textrm{Sing}(T)=
\frac{1}{3}(1,1)+\frac{1}{3}(1,2)$, as required. \\ \\

$\bullet$ Case $(3v)$. $\; G=\mathcal{S}_3 \ltimes
(\mathbb{Z}_4)^2=G(96,64), \; \; \mathbf{m}=(2,3,8),\; \; g(C)=33,
\; \; \textrm{Sing}(T)=\frac{1}{3}(1,1)+\frac{1}{3}(1,2)$. Set
\begin{equation*}
\begin{split}
g_1 &=zxz^3, \; \; g_2=y, \; \; g_3=xyxzxz^3 \\
\ell_1 &=y, \; \; h_1=yz, \; \; h_2=xy.
\end{split}
\end{equation*}
We have $\mathscr{S}=\bigcup_{\sigma \in G}  \langle
\sigma y \sigma^{-1} \rangle \; \cap G^{\times} = \textrm{Cl}(y)$, hence $|\mathscr{S}|=32$.
In fact, $G$ contains precisely
$16$ subgroups of order $3$, which are all conjugate. For all $h \in \mathscr{S}$
\begin{equation*}
\begin{split}
|\textrm{Fix}_{F,1}(h)|&=|\textrm{Fix}_{F,2}(h)|=1 \\
|\textrm{Fix}_{C,1}(h)|&=|\textrm{Fix}_{C,2}(h)|=1.
\end{split}
\end{equation*}
So $C \times F$ contains exactly $64$ points with nontrivial stabilizer
and for each of them the $G$-orbit has
cardinality $|G|/|y|=32$. Looking at the rotation constants we see
that $\textrm{Sing}(T)=
\frac{1}{3}(1,1)+\frac{1}{3}(1,2)$, as required. \\ \\

$\bullet$ Case $(3w)$. $\; G=\textrm{PSL}_2(\mathbb{F}_7), \; \;
\mathbf{m}=(2,3,7),\; \; g(C)=57,
\; \; \textrm{Sing}(T)=\frac{1}{3}(1,1)+\frac{1}{3}(1,2)$. \\
It is well known that $G$ can be embedded in $S_8$; in fact $G= \langle (375)(486), \,
(126)(348) \rangle$. Set
\begin{equation*}
\begin{split}
g_1 &=(12)(34)(58)(67), \; \; g_2=(154)(367), \; \; g_3=(1247358) \\
\ell_1 &=(154)(367), \; \; h_1=(2465837), \; \; h_2=(1352678).
\end{split}
\end{equation*}
We have $\mathscr{S}= \textrm{Cl}((154)(367))$, so $|\mathscr{S}|=56$. In fact, $G$ contains precisely
$28$ subgroups of order $3$, which are all conjugate. For all $h \in \mathscr{S}$
\begin{equation*}
\begin{split}
|\textrm{Fix}_{F,1}(h)|&=|\textrm{Fix}_{F,2}(h)|=1 \\
|\textrm{Fix}_{C,1}(h)|&=|\textrm{Fix}_{C,2}(h)|=1.
\end{split}
\end{equation*}
So $C \times F$ contains exactly $112$ points with nontrivial stabilizer
and for each of them the $G$-orbit has
cardinality $|G|/|(154)(367)|=56$. Looking at the rotation constants
we see that $\textrm{Sing}(T)=
\frac{1}{3}(1,1)+\frac{1}{3}(1,2)$, as required. \\ \\
In all cases $S$ contains only one
 singular Albanese fibre $\bar{F}$, which is illustrated in Figure \ref{fig-K-5} below.
\begin{figure}[ht!]
\begin{center}
\includegraphics*[totalheight=6 cm]{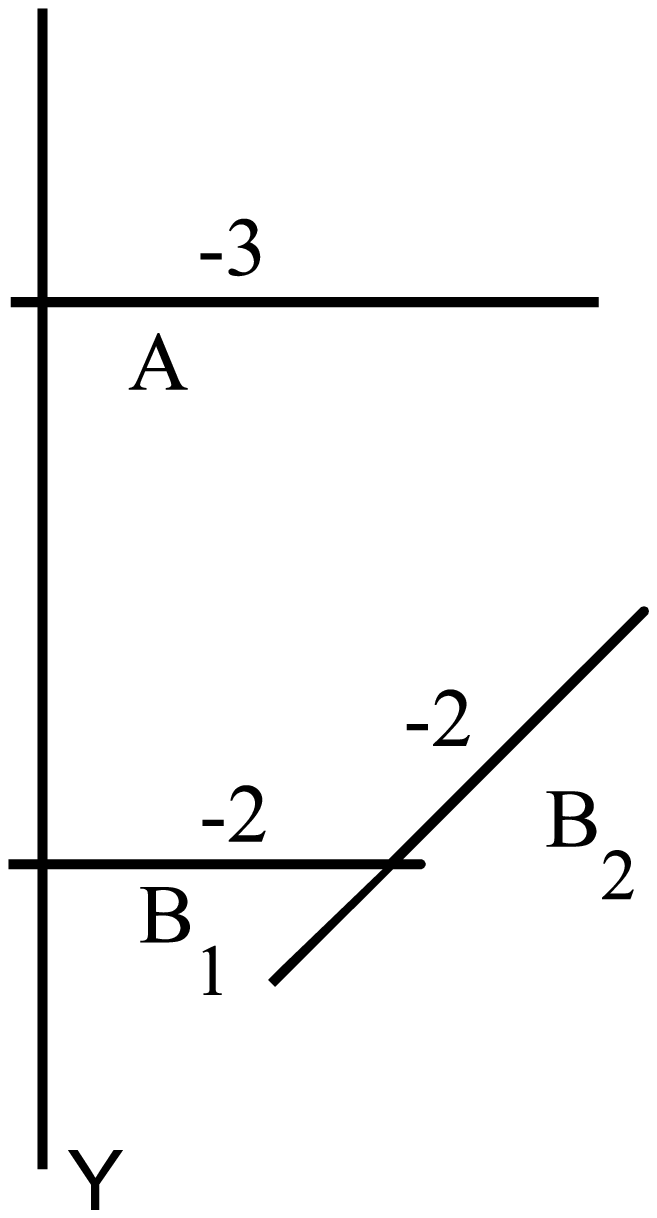}
\end{center}
\caption{}
\label{fig-K-5}
\end{figure}

Here $A$ is a $(-3)$-curve, whereas $B_1$ and $B_2$ are $(-2)$-curves. Since
$\mathbf{n}=(3)$, the central component $Y$ has multiplicity
$3$ in $\bar{F}$ (see Theorem \ref{Serrano}) and a
straightforward computation, using  $\bar{F} A = \bar{F} B_1 =\bar{F} B_2 = 0 $,
shows that
\begin{equation*}
\bar{F}=3Y+A+2B_1+B_2.
\end{equation*}
Using $K_S\bar{F}= 2g(F)-2 = 4$ and $\bar{F}^2=0$ we obtain $K_SY=1$ and $Y^2=-1$.
Hence $Y$ is not a $(-1)$-curve and $S$ is minimal.
\end{proof}

\subsection{The case $K_S^2 = 4$}

\begin{lemma} \label{no-4-gen}
Referring to Table \emph{3} of Appendix B, in cases $(3i)$, $(3j)$, $(3s)$, $(3v)$ the group $G$ is
not $(1 \, | \, 4)$-generated.
\end{lemma}
\begin{proof}
In cases $(3i)$ and $(3j)$ the commutator subgroup $[G,G]$ has order $2$; in case $(3s)$ we have
$[G,G]=\mathcal{A}_4$, which contains no elements of order $4$. In case $(3v)$ we have $G=G(96,64)$;
if $h_1, h_2 \in G$ and $|[h_1,h_2]|=4$ then $|\langle h_1,h_2 \rangle| \leq 48$, so $G$ is not
$(1 \, | \, 4)$-generated.
\end{proof}

\begin{proposition}
Let $\lambda \colon S \lr T=(C \times F)/G$ be a standard
isotrivial fibration with $p_g=q=1$.
If $K_S^2=4$ then $T$ has only \emph{RDPs}.
\end{proposition}
\begin{proof}
Assume that $K_S^2=4$ and $T$ contains at least one singularity
which is not a RDP. Then by Proposition \ref{possibilities}
the only possibility is

\begin{equation*}
g(F)=3, \quad \mathbf{n}=(4), \quad \sT=\frac{1}{2}(1,1)+2 \times
\frac{1}{4}(1,1).
\end{equation*}
In particular $G$ must be one of the groups in Table $3$ of Appendix $B$. Using
Corollary \ref{invariants-S}, we obtain $g(C)-1 = \frac{3}{8}|G|$, so $8$ divides
 $|G|$; moreover,
since $\mathbf{n}=(4)$, it follows that $G$ must be $(1 \, | \,
4)$-generated. Cases $(3i)$, $(3j)$, $(3s)$, $(3v)$ are excluded by
Lemma \ref{no-4-gen}; cases $(3b)$, $(3c)$, $(3g)$, $(3h)$, $(3l)$,
$(3m)$, $(3o)$, $(3p)$ are excluded by Proposition \ref{divides};
cases $(3n)$ and $(3w)$ are excluded because the signature
 $\mathbf{m}$ is not compatible with the singularities of $T$. It remains to rule out cases
  $(3q)$, $(3r)$, $(3t)$, $(3u)$. \\ \\
$\bullet$ Case $(3q)$. $\; G = \mZ_2 \ltimes (\mZ_2 \times
\mZ_8)=G(32,9), \; \; \mathbf{m}=(2,4,8),\; \;
 \textrm{Sing}(T)=\frac{1}{2}(1,1)+ 2 \times \frac{1}{4}(1,1)$. \\
Let $\mathcal{W}=\{\ell_1; \, h_1, \, h_2\}$ be a generating
vector of type $(1 \, | \, 4)$ for $G$. Since $[G, \, G]= \langle yz^2
\rangle$, we may assume $\ell_1=yz^2$. Then $\ell_1 \sim_G \ell_1^{-1}$, hence the same
argument used in proof of Proposition \ref{all-sing-points} shows
 that if $T$ contains a singular point of type $\frac{1}{4}(1,1)$ then it must also
 contain a singular point of type $\frac{1}{4}(1,3)$. Therefore this case cannot occur.
 \\ \\
$\bullet$ Case $(3r)$. $\; G= \mZ_2 \ltimes D_{2,8,5}=G(32,11), \;
\; \mathbf{m}=(2,4,8),\; \;
 \textrm{Sing}(T)=\frac{1}{2}(1,1)+ 2 \times \frac{1}{4}(1,1)$. \\
Let $\mathcal{W}=\{\ell_1; \, h_1, \, h_2\}$. Since we have $[G, \, G]= \langle yz^2
\rangle$, we may assume $\ell_1=yz^2$. Then $\ell_1 \sim_G \ell_1^{-1}$, and
this case can be excluded as the previous one. \\ \\
$\bullet$ Case $(3t)$. $\; G=G(48,33), \; \;
\mathbf{m}=(2,3,12),\; \; \textrm{Sing}(T)=\frac{1}{2}(1,1)+ 2 \times \frac{1}{4}(1,1)$. \\
 We have $[G,G]=Q_8$ and all the elements of order $4$ in $[G,G]$ are conjugate in $G$;
 hence the same argument used in proof of Proposition \ref{all-sing-points} shows
 that if $T$ contains a singular point of type $\frac{1}{4}(1,1)$ then it must also
 contain a singular point of type $\frac{1}{4}(1,3)$. Therefore this case cannot occur.
 \\ \\
$\bullet$ Case $(3u)$. $\; G=\mZ_3 \ltimes (\mZ_4)^2=G(48,3), \; \;
\mathbf{m}=(3^2,4),\; \;
 \textrm{Sing}(T)=\frac{1}{2}(1,1)+ 2 \times \frac{1}{4}(1,1)$. \\
 Let $\mathcal{V}=\{g_1, \, g_2, \, g_3\}$ and $\mathcal{W}=\{\ell_1; \, h_1, \, h_2\}$.
 We have $[G,\,G]=\langle y,z \rangle \cong \mZ_4 \times \mZ_4$ and the conjugacy classes
 in $G$ of elements of order $4$ in $[G, \, G]$ are as follows:
 \begin{equation*}
\begin{split}
\textrm{Cl}(y)&=\{y, \, z, \, y^3z^3 \}, \quad \textrm{Cl}(y^3)=\{y^3, \, z^3, \, yz \} \\
\textrm{Cl}(y^2z)&= \{y^2z, \, y^3z, \, y^3z^2 \}, \quad \textrm{Cl}(yz^2)=\{yz^2, \,
yz^3, \, y^2z^3\}.
\end{split}
\end{equation*}
If  $\ell_1 \sim_G g_3$ then $T$ contains only singularities of type
$\frac{1}{4}(1,1)$, whereas if $\ell_1 \sim_G g_3^{-1}$ then $T$  contains only
singularities of type $\frac{1}{4}(1,3)$. Otherwise $T$ contains only singularities of
type $\frac{1}{2}(1,1)$. Therefore this case cannot occur.
\end{proof}

\subsection{The case $K_S^2 = 3$}

\begin{proposition}
Let $\lambda \colon S \lr T=(C \times F)/G$ be a standard
isotrivial fibration with $p_g=q=1$, $K_S^2=3$ such that
$T$ contains at least one singularity which is not a \emph{RDP}.
Then
\begin{equation*}
g(F)=2, \quad \mathbf{n}=(6), \quad  \ST=2 \times \si{2}{1}
+\si{3}{1} + \si{3}{2}.
\end{equation*}
Furthermore exactly the
following cases occur:
\begin{center}
\begin{tabular}{|c|c|c|c|c|}
\hline
$ $ & \verb|IdSmall| & $ $ & & \\
$G$ & \verb|Group|$(G)$ & $\mathbf{m}$ & $g(C)$ & \emph{Is} $S$ \emph{minimal?}\\
\hline \hline $\mathbb{Z}_2 \ltimes ((\mathbb{Z}_2)^2 \times
\mathbb{Z}_3)  $ & $G(24,8)$ & $(2, 4, 6)$ & $11$ & \emph{Yes} \\
\hline
$\emph{GL}_2(\mathbb{F}_3)$ & $G(48,29)$ & $(2,3,8)$  & $21$ & \emph{Yes} \\
\hline
\end{tabular}
\end{center}
\end{proposition}

\begin{proof}
If $K_S^2=3$ then by Proposition \ref{possibilities} there are two possibilities, namely
\begin{equation*}
\begin{split}
& (a)\; \; g(F)=2, \quad \mathbf{n}=(2,4), \quad
\sT= 2 \times \frac{1}{2}(1,1)+ 2\times \frac{1}{4}(1,1); \\
& (b)\; \; g(F)=2, \quad \mathbf{n}=(6), \quad \quad \sT=2 \times
\frac{1}{2}(1,1)+ \times \frac{1}{3}(1,1)+\frac{1}{3}(1,2).
\end{split}
\end{equation*}
In particular $G$ must be one of the groups in Table $2$ of Appendix $B$.
 We refer to this table and we consider separately the two cases. \\ \\
$\mathbf{Case \;(a).}$ Using Corollary \ref{invariants-S} we obtain $g(C)-1= \frac{5}{8}
|G|$, so $8$ divides $|G|$; moreover, since $\mathbf{n}=(2,4)$, it follows that $|G|$
must be $(1\,|\,2,4)$-generated. Cases $(2b)$, $(2c)$, $(2g)$ are excluded by Proposition
\ref{divides}, whereas case $(2i)$ is excluded because $\textrm{GL}_2(\mathbb{F}_3)$ is
not $(1\,|\,2,4)$-generated (this can be easily checked with \verb|GAP4|).
In cases $(2f)$ and $(2h)$ each element of order $4$ in $G$ is conjugate
to its inverse, hence the same argument used in proof of
Proposition \ref{all-sing-points}
shows that if $T$ contains a singular point of type $\frac{1}{4}(1,1)$ then it must also
 contain a singular point of type $\frac{1}{4}(1,3)$. Therefore this case cannot occur.
 \\ \\
$\mathbf{Case \;(b).}$ Using Corollary \ref{invariants-S} we obtain
$g(C)-1= \frac{5}{12}
|G|$, so $12$ divides $|G|$; moreover, since $\mathbf{n}=(6)$, it follows that $G$ must
be $(1\,|\,6)$-generated. Cases $(2d)$, $(2e)$, $(2h)$ are excluded by Proposition
\ref{divides}; it remains to show that cases $(2g)$ and $(2i)$ actually occur. \\ \\
$\bullet$ Case $(2g)$. $\;G=\mathbb{Z}_2 \ltimes ((\mathbb{Z}_2)^2
\times \mathbb{Z}_3)=G(24,8), \; \; \mathbf{m}=(2,4,6),\; \;
g(C)=11,
\; \; \textrm{Sing}(T)=2 \times \frac{1}{2}(1,1) + \frac{1}{3}(1,1)+\frac{1}{3}(1,2)$. \\
Set
\begin{equation*}
\begin{split}
g_1 &=x, \; \; g_2=zwx, \; \; g_3=yzw \\
\ell_1 &=yw, \; \; h_1=zw, \; \; h_2=x.
\end{split}
\end{equation*}
We have
\begin{equation*}
\begin{split}
\langle \ell_1 \rangle & = \{ 1, yw, w^2, y, w, yw^2 \} \\
\langle g_2 \rangle & = \{ 1, zwx, y, yzwx \} \\
\langle g_3 \rangle & = \{ 1, yzw, w^2, yz, w, yzw^2 \}.
\end{split}
\end{equation*}
One easily checks that
\begin{itemize}
\item[-] the subgroup $\langle \ell_1 \rangle$ is conjugate only to itself;
\item[-] the subgroup $\langle g_3 \rangle$ is conjugate to
\begin{equation*}
\langle zw^2 \rangle= \{ 1, zw^2, w, z, w^2, zw \};
\end{equation*}
\item[-] there are six subgroups of $G$ conjugate to $\langle g_2 \rangle$ and
 different from it; all of them
contain $Z(G)=\langle y \rangle$ as their unique subgroup of order $2$.
\end{itemize}
Therefore $\mathscr{S}= \textrm{Cl}(y) \cup \textrm{Cl}(w)= \{y, \, w, \, w^2\}$.
Moreover
\begin{equation*}
\begin{split}
|\textrm{Fix}_{F}(y)|&=6 \\
|\textrm{Fix}_{C}(y)|&=4 \\
|\textrm{Fix}_{F,1}(w)|&=|\textrm{Fix}_{F,2}(w)|=2 \\
|\textrm{Fix}_{C,1}(w)|&=|\textrm{Fix}_{C,2}(w)|=2.
\end{split}
\end{equation*}
Hence $C \times F$ contains exactly
\begin{itemize}
\item[-] $24$ points having stabilizer of order $|y|=2$ and $G$-orbit of cardinality $12$;
\item[-] $16$ points having stabilizer of order $|w|=3$ and $G$-orbit of cardinality $8$.
\end{itemize}
Looking at the rotation constants we see that $\textrm{Sing}(T)=
 2 \times \frac{1}{2}(1,1)+\frac{1}{3}(1,1)+\frac{1}{3}(1,2)$, as required, so
 this case occurs. \\ \\
$\bullet$ Case $(2i)$. $\;G=\textrm{GL}_2(\mathbb{F}_3), \; \;
\mathbf{m}=(2,3,8), \; \;
g(C)=21,\; \; \textrm{Sing}(T)=2 \times \frac{1}{2}(1,1) + \frac{1}{3}(1,1)+\frac{1}{3}(1,2) $.\\
Set
\begin{equation*}
g_1=\left(
  \begin{array}{cc}
    1 & \;\;1 \\
    0 & -1 \\
  \end{array}
\right) \quad g_2=\left(
  \begin{array}{cc}
    0 & -1 \\
    1 & -1 \\
  \end{array}
\right)\quad g_3=\left(
  \begin{array}{cc}
   -1 & \;\;1 \\
   -1 & -1 \\
  \end{array}
\right)
\end{equation*}
\begin{equation*}
\;\; \ell_1=\left(
  \begin{array}{cc}
    1 & -1 \\
    1 & \; \; 0 \\
  \end{array}
\right) \quad
 h_1=\left(
  \begin{array}{cc}
    -1 & -1 \\
    -1 & \; \; 0 \\
  \end{array}
\right)\quad h_2=\left(
  \begin{array}{cc}
    -1 & \; \; 0 \\
    -1 & -1\\
  \end{array}
\right)
\end{equation*}
and $\ell= \left( \begin{array}{cc}
    -1 & \;\;0 \\
   \;\;0 & -1 \\
  \end{array}
\right)$. We have $(\ell_1)^3=(g_3)^4= \ell$ and $(\ell_1)^2=g_2$. Therefore
$\mathscr{S}=\textrm{Cl}(\ell) \cup \textrm{Cl}(g_2) \cup \textrm{Cl}((g_2)^2)=\{\ell \} \cup \textrm{Cl}(g_2)$.
All the eight elements of order $3$ in $G$ are conjugate, so for all $h \in
\textrm{Cl}(g_2)$  we have
\begin{equation*}
\begin{split}
|\textrm{Fix}_{F,1}(h)|&=|\textrm{Fix}_{F,2}(h)|=2\\
|\textrm{Fix}_{C,1}(h)|&=|\textrm{Fix}_{C,2}(h)|=1.
\end{split}
\end{equation*}
Moreover
\begin{equation*}
|\textrm{Fix}_F(\ell)|=6, \quad |\textrm{Fix}_C(\ell)|=8.
\end{equation*}
Therefore $C \times F$ contains exactly
\begin{itemize}
\item[-] $32$ points having stabilizer of order $|g_2|=3$ and $G$-orbit of cardinality $16$;
\item[-] $48$ points having stabilizer of order $|\ell|=2$ and $G$-orbit of cardinality $24$.
\end{itemize}
Looking at the rotation constants we see that $\textrm{Sing}(T)=
 2 \times \frac{1}{2}(1,1)+\frac{1}{3}(1,1)+\frac{1}{3}(1,2)$, as required. \\ \\
In all cases $S$ contains only one
 singular Albanese fibre $\bar{F}$, which is illustrated in Figure \ref{fig-K-3}.
\begin{figure}[ht!]
\begin{center}
\includegraphics*[totalheight=6 cm]{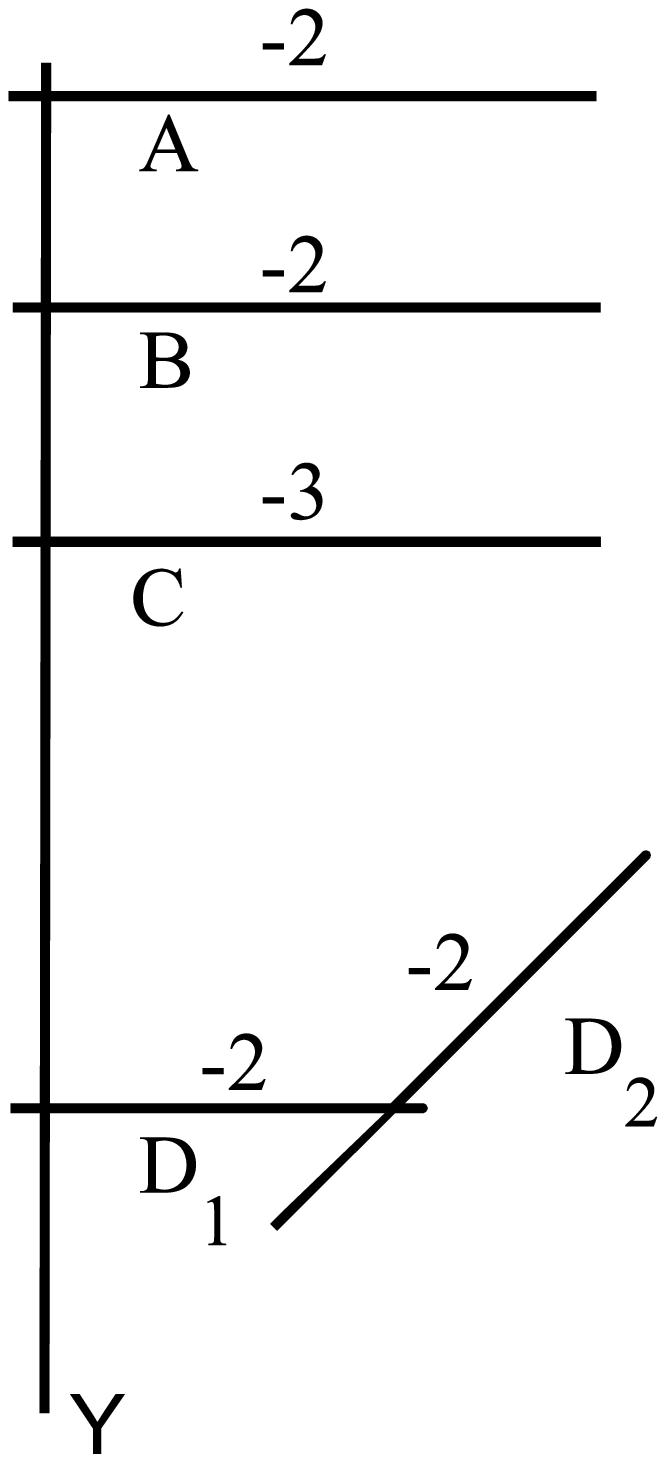}
\end{center}
\caption{} \label{fig-K-3}
\end{figure}
Here $A$, $B$, $D_1$ and $D_2$ are $(-2)$-curves, $C$ is a $(-3)$-curve and a straightforward
computation shows that
\begin{equation*}
\bar{F}=6Y+3A+3B+2C+4D_1+2D_2.
\end{equation*}
Using $K_S\bar{F}=2$ and $\bar{F}^2=0$ we obtain $K_SY=0$ and $Y^2=-2$, so $Y$ is not a
$(-1)$-curve and $S$ is minimal.
\end{proof}

\subsection{The case $K_S^2 = 2$} \label{sec-K2-2}

\begin{proposition}
Let $\lambda \colon S \lr T=(C \times F)/G$ be a standard
isotrivial fibration with $p_g=q=1$, $K_S^2=2$ such that
$T$ contains at least one singularity which is not a \emph{RDP}.
Then there are three possibilities:
\begin{equation*}
\begin{split}
& (d) \; \; g(F)=2, \quad \mathbf{n}=(4),  \quad \ST= 2 \times
\si{2}{1} +\si{4}{1} + \si{4}{3} \\
& (e) \; \; g(F)=3, \quad \mathbf{n}=(4),  \quad \ST= 4 \times
\si{4}{1} \\
& (g) \; \; g(F)=2, \quad \mathbf{n}=(3), \quad \ST= 2 \times
\si{3}{1} + 2 \times \si{3}{2}.
\end{split}
\end{equation*}
In case $(d)$ exactly the following two subcases occur:
\begin{center}
\begin{tabular}{|c|c|c|c|c|}
\hline
$ $ & \verb|IdSmall| & $ $ & &\\
$G$ & \verb|Group|$(G)$ & $\mathbf{m}$ & $g(C)$ &  \emph{Is} $S$ \emph{minimal?} \\
\hline \hline
$D_{2,8,3}$ & $G(16,8)$ & $(2, 4, 8)$ &  $7$ & \emph{Yes} \\
\hline
$\emph{SL}_2(\mathbb{F}_3)$ & $G(24,3)$ & $(3^2,4)$  & $10$ & \emph{Yes} \\
\hline
\end{tabular}
\end{center}
In case $(e)$ there is just one occurrence:
\begin{center}
\begin{tabular}{|c|c|c|c|c|}
\hline
$ $ & \verb|IdSmall| & $ $ & &\\
$G$ & \verb|Group|$(G)$ & $\mathbf{m}$ & $g(C)$ &  \emph{Is} $S$ \emph{minimal?} \\
\hline \hline
$\mZ_3 \ltimes (\mZ_4)^2$ & $G(48,3)$ & $(3^2,4)$ &
$19$ & \emph{No}
\\ \hline
\end{tabular}
\end{center}
Finally, in case $(g)$ there are exactly the subcases below:
\begin{center}
\begin{tabular}{|c|c|c|c|c|}
\hline
$ $ & \verb|IdSmall| & $ $ & &\\
$G$ & \verb|Group|$(G)$ & $\mathbf{m}$ & $g(C)$ &  \emph{Is} $S$ \emph{minimal?} \\
\hline \hline
$\mathcal{S}_3$ & $G(6,1)$ & $(2^2, 3^2)$ & $3$ & \emph{Yes} \\
\hline
$D_{4,3,-1}$ & $G(12,1)$ & $(3, 4^2)$ & $5$ & \emph{Yes} \\
\hline
$D_6$ & $G(12,4)$ & $(2^3, 3)$ & $5$ & \emph{Yes} \\
\hline
\end{tabular}
\end{center}
Moreover in case $(e)$ the minimal model $\widehat{S}$ of $S$
 satisfies $K_{\widehat{S}}^2=3$.
\end{proposition}
\begin{proof}
If $K_S^2=2$ then by Proposition \ref{possibilities} there are seven
possibilities, namely
\begin{equation*}
\begin{split}
& (a)\; \; g(F)=3, \quad \mathbf{n}=(16), \quad \sT=
\frac{1}{2}(1,1)+ \frac{1}{8}(1,1)+ \frac{1}{8}(1,3); \\
& (b)\; \; g(F)=2, \quad \mathbf{n}=(8), \quad \; \; \sT=
\frac{1}{2}(1,1)+ \frac{1}{8}(1,3)+ \frac{1}{8}(1,5); \\
& (c)\; \; g(F)=3, \quad \mathbf{n}=(12), \quad
\sT=\frac{1}{3}(1,2)+ 2 \times \frac{1}{6}(1,1); \\
& (d)\; \; g(F)=2, \quad \mathbf{n}=(4), \quad \;\sT= 2 \times
\frac{1}{2}(1,1)+ \frac{1}{4}(1,1)+ \frac{1}{4}(1,3); \\
& (e)\; \; g(F)=3, \quad \mathbf{n}=(4), \quad \;\sT= 4 \times
\frac{1}{4}(1,1); \\
& (f)\; \; g(F)=2, \; \;  \mathbf{n}=(4^2), \quad \sT= 4 \times
\frac{1}{4}(1,1); \\
& (g) \; \; g(F)=2, \; \; \mathbf{n}=(3), \quad \; \; \sT=2 \times \frac{1}{3}(1,1)
+ 2 \times \frac{1}{3}(1,2).
\end{split}
\end{equation*}
If $g(F)=2$ then $G$ must be one of the groups in Table $2$ of
Appendix $B$, whereas if $g(F)=3$ then $G$ must be one of the groups
in Table $3$. Let us consider separately the different cases. \\ \\
$\mathbf{Case \;(a).}$ Using Corollary \ref{invariants-S} we obtain
$g(C)-1=\frac{15}{32} |G|$, hence $32$ divides $|G|$; looking at
Table $3$ we see that the only possibilities  are
$(3q)$ and $(3r)$. In both cases $[G, \, G]$ has order $4$, so $G$
is not $(1 \, | \, 16)$-generated and this contradicts
$\mathbf{n}=(16)$. Hence this case does not occur. \\ \\
$\mathbf{Case \;(b).}$ We obtain $g(C)-1=\frac{7}{16} |G|$, hence
$16$ divides $|G|$; looking at
 Table $2$ we see that the only possibilities are
  $(2f)$ and $(2i)$. In both cases one
 easily checks that $[G, \, G]$ contains no elements of order
 $8$, so $G$ is not $(1 \, | \, 8)$-generated and this contradicts
$\mathbf{n}=(8)$. Hence this case does not occur. \\ \\
$\mathbf{Case \;(c).}$ We obtain $g(C)-1= \frac{11}{24} |G|$, so
$24$ divides $|G|$. Referring to Table $3$ of Appendix $B$, we are
left with cases $(3l)$, $(3m)$, $(3n)$, $(3o)$, $(3p)$, $(3s)$,
$(3t)$, $(3u)$, $(3v)$, $(3w)$. All these possibilities can be ruled
out by using Proposition \ref{divides}, hence this case does not
occur. \\ \\
$\mathbf{Case \;(d).}$ We obtain $g(C)-1= \frac{3}{8} |G|$, so $8$
divides $|G|$. By direct
computation or using \verb|GAP4| one checks that the groups in
cases $(2b)$, $(2c)$, $(2g)$ and $(2i)$ are not $(1\, | \, 4)$-generated,
contradicting $\mathbf{n}=(4)$; so
the only possibilities are $(2f)$ and
$(2h)$. Let us show that both actually occur. \\ \\
$\bullet$ Case $(2f)$. $\;G=D_{2,8,3}, \; \; \mathbf{m}=(2,4,8),\;
\; g(C)=7,
\; \; \textrm{Sing}(T)=2 \times \frac{1}{2}(1,1) + \frac{1}{4}(1,1)+\frac{1}{4}(1,3)$. \\
Set
\begin{equation*}
\begin{split}
g_1 &=x, \; \; g_2=xy^7, \; \; g_3=y \\
\ell_1 &=y^2, \; \; h_1=y, \; \; h_2=x.
\end{split}
\end{equation*}
We have
\begin{equation*}
\textrm{Cl}(y)=\{y, \, y^3\}, \; \; \; \textrm{Cl}(y^2)=\{y^2, \,
y^6\}, \; \; \; \textrm{Cl}(y^4)=\{y^4\}.
\end{equation*}
Since $(g_2)^2=(\ell_1)^2$ and $(g_3)^2=\ell_1$, we obtain
$\mathscr{S}=\bigcup_{\sigma \in G} \langle \sigma y^2 \sigma ^{-1} \rangle
 \;\cap G ^{\times}= \{y^2, \, y^4,
\, y^6\}$. Moreover
\begin{equation*}
\begin{split}
|\textrm{Fix}_{F}(y^4)|&=6 \\
|\textrm{Fix}_{C}(y^4)|&=4 \\
|\textrm{Fix}_{F,1}(y^2)|&=|\textrm{Fix}_{F,3}(y^2)|=1 \\
|\textrm{Fix}_{C,1}(y^2)|&=|\textrm{Fix}_{C,3}(y^2)|=2.
\end{split}
\end{equation*}
Therefore $C \times F$ contains exactly
\begin{itemize}
\item[-] $16$ points having stabilizer of order $|y^4|=2$ and $G$-orbit of cardinality $8$;
\item[-] $8$ points having stabilizer of order $|y^2|=4$ and $G$-orbit of cardinality $4$.
\end{itemize}
Looking at the rotation constants we see that $\textrm{Sing}(T)=
 2 \times \frac{1}{2}(1,1)+\frac{1}{4}(1,1)+\frac{1}{4}(1,3)$, as required.
\\ \\
$\bullet$ Case $(2h)$. $\;G=\textrm{SL}_2(\mathbb{F}_3), \; \;
\mathbf{m}=(3^2,4), \; \;
g(C)=10,\; \; \textrm{Sing}(T)=2 \times \frac{1}{2}(1,1) + \frac{1}{4}(1,1)+\frac{1}{4}(1,3) $.\\
Set
\begin{equation*}
g_1=\left(
  \begin{array}{cc}
    \;\;\;0 & \;\;1 \\
    -1 & -1 \\
  \end{array}
\right) \quad g_2=\left(
  \begin{array}{cc}
    0 & -1 \\
    1 & -1 \\
  \end{array}
\right)\quad g_3=\left(
  \begin{array}{cc}
   -1 & 1 \\
    \;\;1 & 1 \\
  \end{array}
\right)
\end{equation*}
\begin{equation*}
\;\; \ell_1=\left(
  \begin{array}{cc}
    -1 & 1 \\
    \;\;1 & 1 \\
  \end{array}
\right) \quad
 h_1=\left(
  \begin{array}{cc}
    \;\;\;0 & 1 \\
    -1 & 0 \\
  \end{array}
\right)\quad h_2=\left(
  \begin{array}{cc}
    1 & 1 \\
    0 & 1\\
  \end{array}
\right)
\end{equation*}
and $\ell= \left( \begin{array}{cc}
    -1 & \;\;0 \\
   \;\;0 & -1 \\
  \end{array}
\right)$. The group $G$ contains six elements of order $4$, which
are all conjugate. Therefore there are three cyclic subgroups $H_1,
\, H_2, \, H_3$ of order $4$, all conjugate and such that $H_i \cap
H_j=\langle \ell \rangle$ for $i \neq j$. If $h \in G $ and $|h|=4$ then
\begin{equation*}
\begin{split}
|\textrm{Fix}_{F,1}(h)|&=|\textrm{Fix}_{F,3}(h)|=1 \\
|\textrm{Fix}_{C,1}(h)|&=|\textrm{Fix}_{C,3}(h)|=1.
\end{split}
\end{equation*}
Therefore $C \times F$ contains exactly
\begin{itemize}
\item[-] $24$ points having stabilizer of order $|\ell|=2$ and $G$-orbit of cardinality $12$;
\item[-] $12$ points having stabilizer of order $|h|=4$ and $G$-orbit of cardinality $6$.
\end{itemize}
Looking at the rotation constants we see that $\textrm{Sing}(T)=
 2 \times \frac{1}{2}(1,1)+\frac{1}{4}(1,1)+\frac{1}{4}(1,3)$, as
 required. \\ \\
Now we show that all surfaces in $\mathbf{Case \; (d)}$ are minimal.
In fact they contain only one
 singular Albanese fibre $\bar{F}$, which is illustrated in Figure \ref{fig-K-2a}.
\begin{figure}[ht!]
\begin{center}
\includegraphics*[totalheight=6 cm]{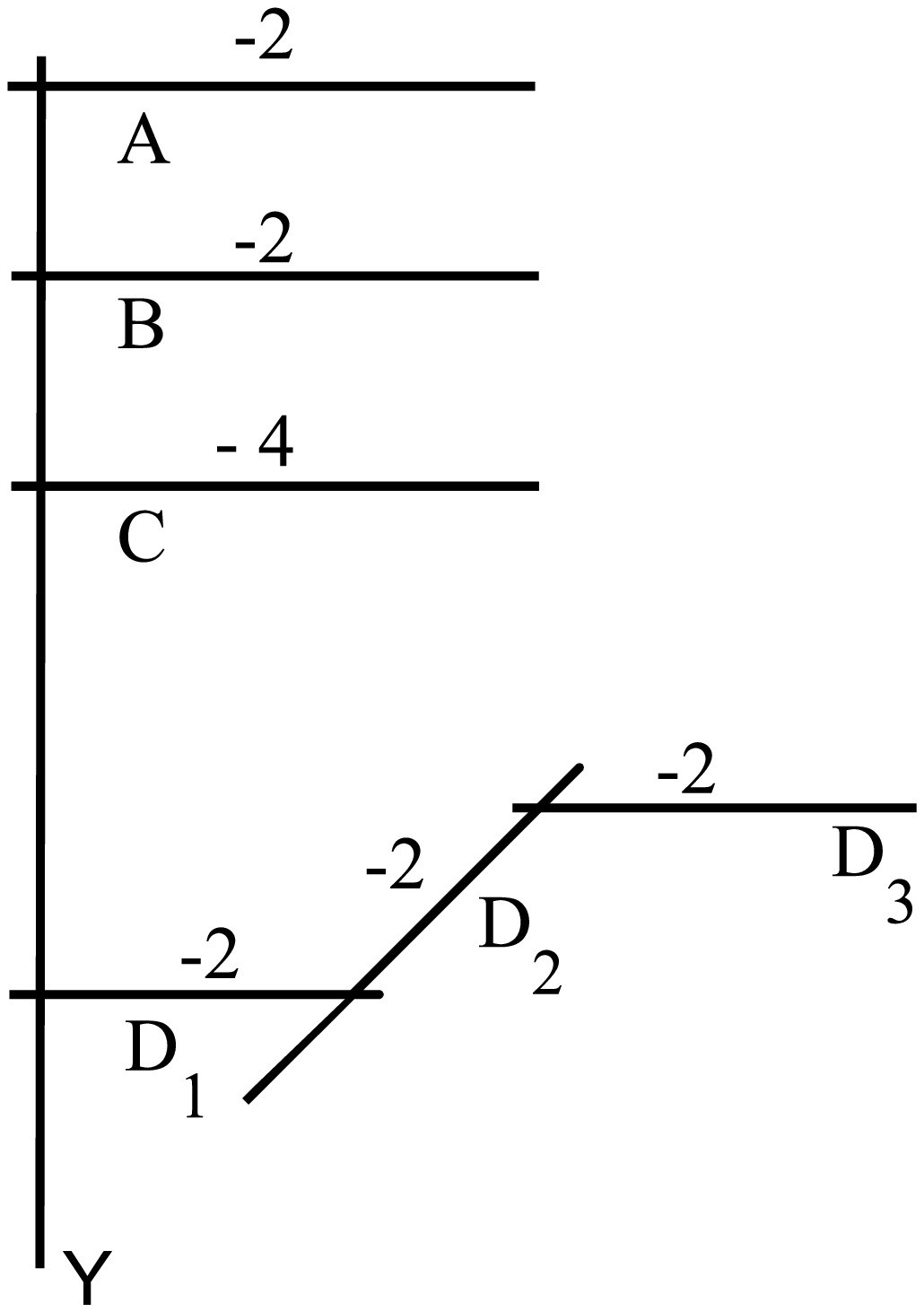}
\end{center}
\caption{} \label{fig-K-2a}
\end{figure}
Here $A$, $B$, $D_1$, $D_2$ and $D_3$ are $(-2)$-curves, $C$ is a
$(-4)$-curve and a straightforward
computation shows that
\begin{equation*}
\bar{F}=4Y+2A+2B+C+3D_1+2D_2+D_3.
\end{equation*}
Using $K_S\bar{F}=2$ and $\bar{F}^2=0$ we obtain $K_SY=0$ and $Y^2=-2$, so $Y$ is not a
$(-1)$-curve and $S$ is minimal. \\ \\
$\mathbf{Case \;(e).}$ We obtain $g(C)-1=\frac{3}{8}|G|$, hence $8$
divides $|G|$. Referring to Table $3$ in Appendix $B$, we have what
follows. \\ \\
- Cases $(3b)$, $(3c)$, $(3g)$, $(3h)$, $(3i)$, $(3j)$, $(3l)$,
$(3m)$,  $(3s)$, $(3v)$ must be excluded because the corresponding
$G$ are not $(1 \, | \, 4)$-generated, contradicting
$\mathbf{n}=(4)$. \\ \\
- Cases $(3n)$ and $(3w)$ must be excluded
because no component of $\mathbf{m}$ is divided by $4$, a
contradiction because the
singularities of $T$ must be of type $\frac{1}{4}(1,1)$. \\ \\
- In cases $(3o)$, $(3p)$, $(3q)$, $(3r)$, $(3t)$ all elements of
order $4$ in $[G,\,G]$ are conjugate in $G$; therefore the same
argument used in proof of Proposition \ref{all-sing-points} shows
that $S$ must contain both $\frac{1}{4}(1,1)$ and $\frac{1}{4}(1,3)$
singularities, a contradiction. \\ \\
Now we show that Case $(3u)$ occurs. \\ \\
$\bullet$ Case $(3u)$. $\;G=\mZ_3 \ltimes (\mZ_4)^2=G(48,3), \; \;
\mathbf{m}=(3^2,4),\; \; g(C)=19,
\; \; \textrm{Sing}(T)=4 \times \frac{1}{4}(1,1)$. \\
Set
\begin{equation*}
\begin{split}
g_1 &=x, \; \; g_2=x^2y^3, \; \; g_3=y \\
\ell_1 &=y, \; \; h_1=x, \; \; h_2=xyxy^2.
\end{split}
\end{equation*}
We have
$\mathscr{S}=\bigcup_{\sigma \in G} \langle \sigma y \sigma ^{-1} \rangle$ and
the elements of order $4$ in $\mathscr{S}$ are precisely
$\{y, \; z, \; y^3z^3, \; y^3, \; z^3, \; yz \}$.
Moreover $\textrm{Cl}(y)=\{y, \; z, \; y^3z^3 \}$ and
$\textrm{Cl}(y^3)=\{y^3, \; z^3, \; yz \}$.
Take any $h \in \mathscr{S}$ such that $|h|=4$; since $h$ is not
conjugate to $h^{-1}$ in $G$, Proposition \ref{fixed-points}
implies
\begin{equation*}
\begin{split}
|\textrm{Fix}_{F,1}(h)|&=4, \quad |\textrm{Fix}_{F,3}(h)|=0, \\
|\textrm{Fix}_{C,1}(h)|&=4, \quad |\textrm{Fix}_{C,3}(h)|=0.
\end{split}
\end{equation*}
Therefore $C \times F$ contains exactly $48$ points with  nontrivial stabilizer
and for each of them the $G$-orbit has cardinality $|G|/|y|=12$.
Looking at the rotation constants we see that $\textrm{Sing}(T)=
 4 \times \frac{1}{4}(1,1)$, as required. The surface $S$ contain only one
 singular Albanese fibre $\bar{F}$, which is illustrated in Figure \ref{fig-K-2}.
\begin{figure}[ht!]
\begin{center}
\includegraphics*[totalheight=6 cm]{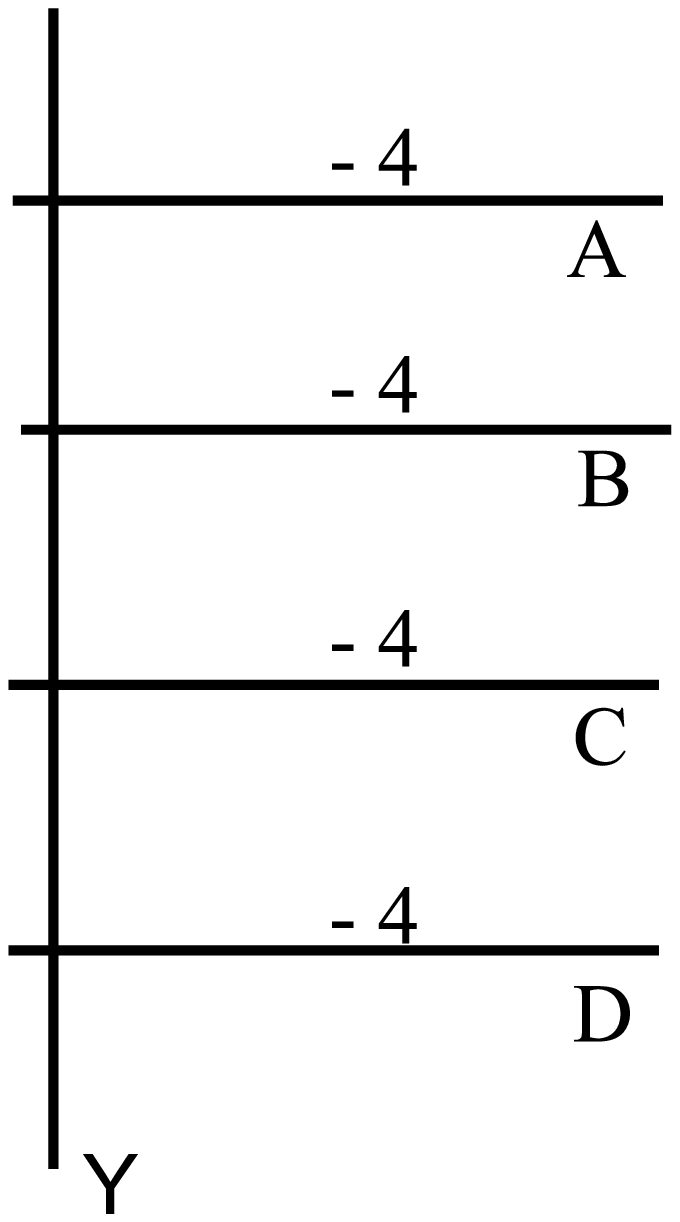}
\end{center}
\caption{} \label{fig-K-2}
\end{figure}
Here $A$, $B$, $C$, $D$ are $(-4)$-curves and a straightforward computation shows
that
\begin{equation*}
\bar{F}=4Y+A+B+C+D.
\end{equation*}
Since $K_S\bar{F}=4$ and $(\bar{F})^2=0$ we obtain $K_SY=Y^2=-1$, i. e.
 $Y$ is the unique $(-1)$-curve in $S$. The minimal model $\widehat{S}$ of $S$ is obtained
by contracting $Y$, hence $K_{\widehat{S}}^2=3$. Therefore $\widehat{S}$ is an example
of minimal surface of general type with $p_g=q=1$, $K^2=g_{\textrm{alb}}=3$ and
a unique singular Albanese fibre. The existence of such surfaces
 was previously established, in a completely different way, by Ishida in
\cite{Is06}. \\ \\
$\mathbf{Case \;(f).}$ We obtain $g(C)-1=\frac{3}{4}|G|$, hence $4$
divides $|G|$; moreover
 $G$ must be $(1\, | \, 4^2)$-generated. Look at Table $2$ of
Appendix $B$. Cases $(2b)$ and $(2e)$ are excluded by using
Proposition \ref{divides}, whereas Case $(2i)$ is excluded because
$\textrm{GL}_2(\mathbb{F}_3)$ is not $(1\, | \, 4^2)$-generated.
In cases
$(2c)$, $(2f)$, $(2g)$ and $(2h)$ all the elements of order $4$ in $G$ are conjugate
to their inverse, hence if $S$ contains a singular point of
type $\frac{1}{4}(1,1)$ it should also contain a singular point
of type $\frac{1}{4}(1,3)$, a contradiction. Hence we must only consider $(2d)$.
In this
case $G=D_{4,3,-1}$, which contains two conjugacy classes of
elements of order $4$, namely $\textrm{Cl}(x)=\{x, \, xy, \, xy^2
\}$ and $\textrm{Cl}(x^3)= \{x^3, \, x^3y, \, x^3y^2 \}$. Since
the only element of order $2$ in $G$ is $x^2$, two
different $2$-Sylow of $G$ intersect exactly in $\langle x^2
\rangle$. This show that $T$ should contain some singular points of
type $\frac{1}{2}(1,1)$, a contradiction.
\\ \\
$\mathbf{Case \;(g).}$ We obtain $g(C)-1=\frac{1}{3}|G|$, hence $3$
divides $|G|$; moreover
 $G$ must be $(1\, | \, 3)$-generated. Referring to Table $2$ of
Appendix $B$, the groups in cases $(2g)$, $(2h)$, $(2i)$ are excluded because
 they are not
$(1\, | \, 3)$-generated, so we are
 left to show that cases $(2a)$, $(2c)$ and $(2e)$ occur. \\ \\
$\bullet$ Case $(2a)$. $\;G=\mathcal{S}_3, \; \;
\mathbf{m}=(2^2,3^2),\; \; g(C)=3,
\; \; \textrm{Sing}(T)=2 \times \frac{1}{3}(1,1) + 2 \times \frac{1}{3}(1,2)$. \\
Set
\begin{equation*}
\begin{split}
g_1 &=(12), \; \; g_2=(12), \; \; g_3=(123), \; \; g_4=(132) \\
\ell_1 &=(123), \; \; h_1=(13), \; \; h_2=(12).
\end{split}
\end{equation*}
We have $\mathscr{S}= \textrm{Cl}((123))=\{(123), \; (132) \}$
and for all $h \in \mathscr{S}$
\begin{equation*}
\begin{split}
|\textrm{Fix}_{F,1}(h)|&= |\textrm{Fix}_{F,2}(h)|=2 \\
|\textrm{Fix}_{C,1}(h)|&= |\textrm{Fix}_{C,2}(h)|=1.
\end{split}
\end{equation*}
Hence $C \times F$ contains exactly $8$ points with nontrivial stabilizer
and the $G$-orbit of each of these points has cardinality $|G|/|(123)|=2$.
Looking at the rotation constants we see that
$\sT=2 \times \frac{1}{3}(1,1) + 2 \times \frac{1}{3}(1,2)$, as required.
\\ \\
$\bullet$ Case $(2d)$. $\;G=D_{4,3,-1}, \; \;
\mathbf{m}=(3,4^2),\; \; g(C)=5,
\; \; \textrm{Sing}(T)=2 \times \frac{1}{3}(1,1) + 2 \times \frac{1}{3}(1,2)$. \\
Set
\begin{equation*}
\begin{split}
g_1 &=y, \; \; g_2=y^2x^3, \; \; g_3=x \\
\ell_1 &=y, \; \; h_1=y, \; \; h_2=x.
\end{split}
\end{equation*}
We have $\mathscr{S}=\textrm{Cl}(y)=\{y, \; y^2 \}$
and for all $h \in \mathscr{S}$
\begin{equation*}
\begin{split}
|\textrm{Fix}_{F,1}(h)|&= |\textrm{Fix}_{F,2}(h)|=2 \\
|\textrm{Fix}_{C,1}(h)|&= |\textrm{Fix}_{C,2}(h)|=2.
\end{split}
\end{equation*}
Hence $C \times F$ contains exactly $16$ points with nontrivial stabilizer
and the $G$-orbit of each of these points has cardinality $|G|/|y|=4$.
Looking at the rotation constants we see that
$\sT=2 \times \frac{1}{3}(1,1) + 2 \times \frac{1}{3}(1,2)$, as required.
\\ \\
$\bullet$ Case $(2e)$. $\;G=D_6, \; \;
\mathbf{m}=(2^3,3),\; \; g(C)=5,
\; \; \textrm{Sing}(T)=2 \times \frac{1}{3}(1,1) + 2 \times \frac{1}{3}(1,2)$. \\
Set
\begin{equation*}
\begin{split}
g_1 &=x, \; \; g_2=xy, \; \; g_3=y^3, \; \; g_4=y^2 \\
\ell_1 &=y^2, \; \; h_1=x, \; \; h_2=y.
\end{split}
\end{equation*}
We have $\mathscr{S}=\textrm{Cl}(y^2)=\{y^2, \; y^4 \}$
and for all $h \in \mathscr{S}$
\begin{equation*}
\begin{split}
|\textrm{Fix}_{F,1}(h)|&= |\textrm{Fix}_{F,2}(h)|=2, \\
|\textrm{Fix}_{C,1}(h)|&= |\textrm{Fix}_{C,2}(h)|=2.
\end{split}
\end{equation*}
Hence $C \times F$ contains exactly $16$ points with nontrivial stabilizer
and the $G$-orbit of each of these points has cardinality $|G|/|y^2|=4$.
Looking at the rotation constants we see that
$\sT=2 \times \frac{1}{3}(1,1) + 2 \times \frac{1}{3}(1,2)$, as required.
\\ \\
Now we show that the surfaces in $\mathbf{Case \;(g)}$ are minimal.
In fact they all contain only one
 singular Albanese fibre $\bar{F}$, which is illustrated in Figure \ref{fig-K-2g}.
\begin{figure}[ht!]
\begin{center}
\includegraphics*[totalheight=6 cm]{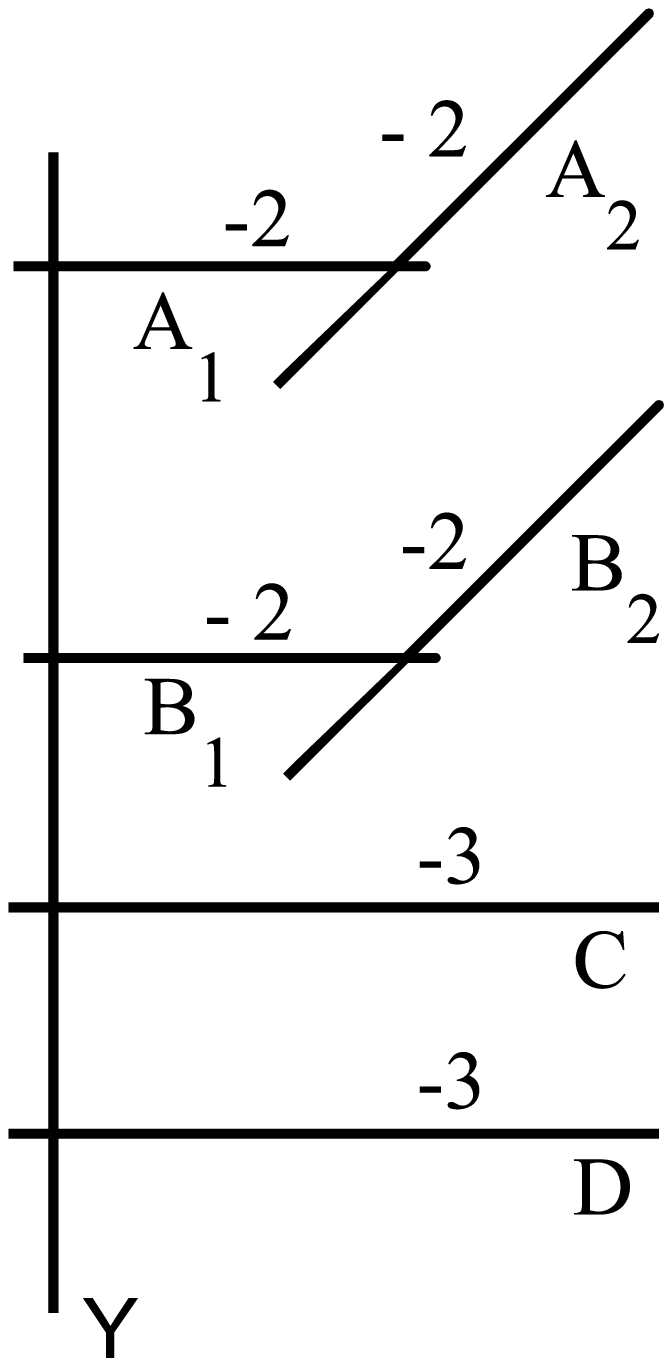}
\end{center}
\caption{} \label{fig-K-2g}
\end{figure}
Here $A_1$, $A_2$, $B_1$, $B_2$ are $(-2)$-curves, $C$, $D$ are
$(-3)$-curves and a straightforward
computation shows that
\begin{equation*}
\bar{F}=3Y+2A_1+A_2+2B_1+B_2+C+D.
\end{equation*}
Using $K_S\bar{F}=2$ and $\bar{F}^2=0$ we obtain $K_SY=0$ and $Y^2=-2$, so $Y$ is not a
$(-1)$-curve and $S$ is minimal.
\end{proof}

\newpage

\section{The case where $S$ is not minimal}

The description of all non minimal examples would put an end to the classification
of standard isotrivial fibrations with $p_g=q=1$; however, it seems to us difficult to
achieve it by using our methods. We can prove the following
\begin{proposition} \label{-1-curves}
Let $\lambda \colon S \lr (C \times F)/G$ be a standard isotrivial fibration
of general type
with $p_g=q=1$. Then $S$ contains at most five $(-1)$-curves.
\end{proposition}
\begin{proof}
Let $\alpha \colon S \lr E$ be the Albanese map of $S$,
let $\widehat{S}$ be the minimal model of $S$ and
$\hat{\alpha} \colon \widehat{S} \lr E$ the Albanese map of $\widehat{S}$.
By Theorem \ref{Serrano} the $(-1)$-curves of $S$ may only appear as central components
of reducible fibres of $\alpha$. Therefore the number of such curves is smaller than or equal to
the number of singular fibres of $\hat{\alpha}$. On the other hand,
 by Zeuthen-Segre formula (\cite[p.116]{Be2}) we have
 \begin{equation*}
10 \geq e(\widehat{S}) = \sum_{x \in \textrm{Crit}(\hat{\alpha})} \mu_x,
\end{equation*}
where $\textrm{Crit}(\hat{\alpha})$ is the set of points of $E$ where the fibre of $\hat{\alpha}$ is
singular. The integer $\mu_x$ satisfies $\mu_x \geq 1$ and
 equality holds if and only if the fibre of $\hat{\alpha}$ over $x$ has
an ordinary node as unique singularity. This would imply that
the general fibre of $\alpha$ is rational, a contradiction. Therefore
$\mu_x \geq 2$ for every $x \in \textrm{Crit}(\hat{\alpha})$, so $\hat{\alpha}$ has at most
five singular fibres.
\end{proof}
The main problem is that further $(-1)$-curves may appear after contracting
the $(-1)$-curves of $S$. This happens for instance in the following example.

\subsection{An example with $K_S^2=1$} \label{non-min-K2-1}

In this section we construct a standard isotrivial fibration $S$
  with $p_g=q=1$ and $K_S^2=1$, whose minimal model $\widehat{S}$ satisfies
$K_{\widehat{S}}^2=3$. The
building data for $S$ are
\begin{equation*}
\begin{split}
g(F)& =3, \quad \mathbf{m}=(3^2,7), \\
g(C)& =10, \quad \mathbf{n}=(7), \\
G & =D_{3,7,2}= \langle x,y \; | \; x^3=y^7=1,\, xyx^{-1}=y^2 \rangle.
\end{split}
\end{equation*}
Set
\begin{equation*}
\begin{split}
g_1 &=x^2, \; \; g_2=xy^6, \; \; g_3=y \\
\ell_1 &=y, \; \; h_1=y, \; \; h_2=x.
\end{split}
\end{equation*}
We have
$\mathscr{S}=\bigcup_{\sigma \in G} \langle \sigma y \sigma ^{-1} \rangle =
\{y, \, y^2, \, y^3, \, y^4, \, y^5, \, y^6 \}$ and moreover
$\textrm{Cl}(y)=\{y, \, y^2, \, y^4 \}$, $\textrm{Cl}(y^3)=\{y^3, \, y^6, \, y^5 \}$. Hence for all $h \in \mathscr{S}$ we obtain
\begin{equation*}
\begin{split}
|\textrm{Fix}_{F,1}(h)|&=|\textrm{Fix}_{F,2}(h)|=|\textrm{Fix}_{F,4}(h)|=1\\
|\textrm{Fix}_{F,3}(h)|&=|\textrm{Fix}_{F,5}(h)|=|\textrm{Fix}_{F,6}(h)|=0\\
|\textrm{Fix}_{C,1}(h)|&=|\textrm{Fix}_{C,2}(h)|=|\textrm{Fix}_{C,4}(h)|=1 \\
|\textrm{Fix}_{C,3}(h)|&=|\textrm{Fix}_{C,5}(h)|=|\textrm{Fix}_{C,6}(h)|=0.
\end{split}
\end{equation*}
It follows that $C \times F$ contains exactly $9$ points
with nontrivial stabilizer and for each of them
the $G$-orbit has cardinality $|G|/|y|=3$. Looking
 at the rotation constants we see that
\begin{equation*}
\textrm{Sing}(T)= \frac{1}{7}(1,1)+\frac{1}{7}(1,2)+ \frac{1}{7}(1,4),
\end{equation*}
so using Proposition \ref{strukture} one checks that
$S$ is a surface of general type with $p_g=q=1$, $K_S^2=1$. Furthermore the surface $S$ contains only one singular Albanese fibre $\bar{F}$ which is
illustrated in Figure \ref{fig-K-1} below.
\begin{figure}[H]
\begin{center}
\includegraphics*[totalheight=6 cm]{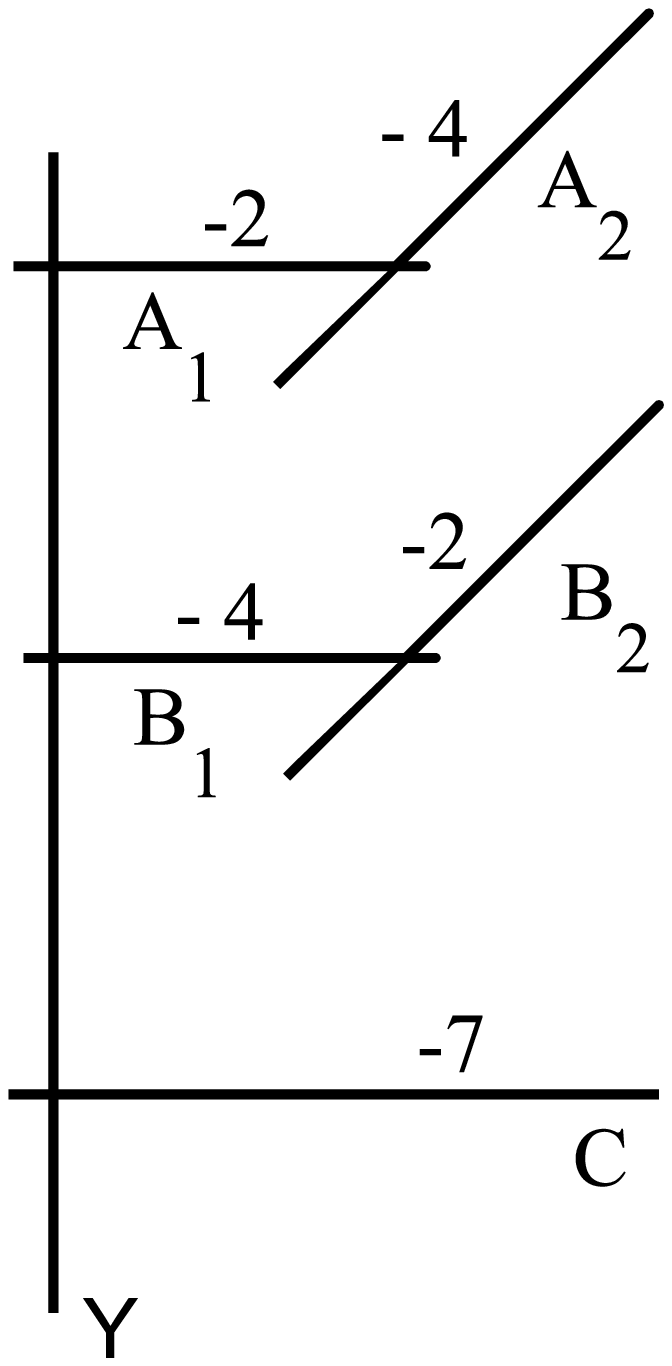}
\caption{} \label{fig-K-1}
\end{center}
\end{figure}
Notice that, since $2 \cdot 4 \equiv 1 \; (\textrm{mod}\; 7)$,
 the cyclic quotient singularities $\frac{1}{7}(1,2)$ and $\frac{1}{7}(1,4)$ are
 analytically isomorphic (see Section \ref{HJ-res}); moreover,
 the resolution algorithm given in
\cite[Chapter II]{Lau71} implies that the
corresponding Hirzebruch-Jung strings are
  attached in a mirror-like way to the central component $Y$ of $\bar{F}$.
A straightforward computation shows that
\begin{equation*}
\bar{F}= 7Y+4 A_1+A_2 +2B_1+B_2 +C.
\end{equation*}
Using $K_S\bar{F}=4$ and $\bar{F}^2=0$ we obtain $K_SY=Y^2=-1$, hence $Y$ is
 the unique $(-1)$-curve in $S$.
The minimal model $\widehat{S}$ of $S$ is obtained by first contracting $Y$ and then the
 image of $A_1$; therefore $K_{\widehat{S}}^2=3$.

\newpage

\section*{Appendix $A$}
List of all cyclic quotient singularities $x=\frac{1}{n}(1,q)$ with
$3 \leq B_x \leq 12$.
\label{tabella}
$$
\begin{tabular}{|c|c|c|c|c|}

\ssi{n}{q}{n/q = [b_1, \dots , b_s]}{q^{\prime}}{B_{\si{n}{q}}}{h_{\si{n}{q}}}

 & & & &  \\

\ssi{2}{1}{[ 2 ]}{1}{3+0}{0}

\ssi{3}{1}{[ 3 ]}{1}{3+2/3}{-1/3}

\ssi{3}{2}{[ 2, 2 ]}{2}{5+1/3}{0}

\ssi{4}{1}{[ 4 ]}{1}{4+1/2}{-1}

\ssi{4}{3}{[ 2, 2, 2 ]}{3}{7+1/2}{0}

\ssi{5}{1}{[ 5 ]}{1}{5+2/5}{-9/5}

\ssi{5}{2}{[ 3, 2 ]}{3}{6+0}{-2/5}

\ssi{5}{4}{[ 2, 2, 2, 2 ]}{4}{9+3/5}{0}

\ssi{6}{1}{[ 6 ]}{1}{6+1/3}{-8/3}

\ssi{6}{5}{[ 2, 2, 2, 2, 2 ]}{5}{11+2/3}{0}

\ssi{7}{1}{[ 7 ]}{1}{7+2/7}{-25/7}

\ssi{7}{2}{[ 4, 2 ]}{4}{6+6/7}{-8/7}

\ssi{7}{3}{[ 3, 2, 2 ]}{5}{8+1/7}{-3/7}

\ssi{8}{1}{[ 8 ]}{1}{8+1/4}{-9/2}

\ssi{8}{3}{[ 3, 3 ]}{3}{6+3/4}{-1}

\ssi{8}{5}{[ 2, 3, 2 ]}{5}{8+1/4}{-1/2}

\ssi{9}{1}{[ 9 ]}{1}{9+2/9}{-49/9}

\ssi{9}{2}{[ 5, 2 ]}{5}{7+7/9}{-2}

\ssi{9}{4}{[ 3, 2, 2, 2 ]}{7}{10+2/9}{-4/9}

\ssi{10}{1}{[ 10 ]}{1}{10+1/5}{-32/5}

\ssi{10}{3}{[ 4, 2, 2 ]}{7}{9+0}{-6/5}

\ssi{11}{1}{[ 11 ]}{1}{11+2/11}{-81/11}

\ssi{11}{2}{[ 6, 2 ]}{6}{8+8/11}{-32/11}

\ssi{11}{3}{[ 4, 3 ]}{4}{7+7/11}{-20/11}

\ssi{11}{7}{[ 2, 3, 2, 2 ]}{8}{10+4/11}{-6/11}

\ssi{12}{5}{[ 3, 2, 3 ]}{5}{8+5/6}{-1}

\ssi{12}{7}{[ 2, 4, 2 ]}{7}{9+1/6}{-4/3}

\ssi{13}{2}{[ 7, 2 ]}{7}{9+9/13}{-50/13}

\ssi{13}{3}{[ 5, 2, 2 ]}{9}{9+12/13}{-27/13}

\ssi{13}{4}{[ 4, 2, 2, 2 ]}{10}{11+1/13}{-16/13}

\ssi{13}{5}{[ 3, 3, 2 ]}{8}{9+0}{-15/13}

\ssi{14}{3}{[ 5, 3 ]}{5}{8+4/7}{-19/7}

\ssi{15}{2}{[ 8, 2 ]}{8}{10+2/3}{-24/5}

\ssi{15}{4}{[ 4, 4 ]}{4}{8+8/15}{-8/3}

\ssi{16}{3}{[ 6, 2, 2 ]}{11}{10+7/8}{-3}

\ssi{16}{7}{[ 3, 2, 2, 3 ]}{7}{10+7/8}{-1}

\ssi{16}{9}{[ 2, 5, 2 ]}{9}{10+1/8}{-9/4}

\ssi{17}{2}{[ 9, 2 ]}{9}{11+11/17}{-98/17}

\ssi{17}{3}{[ 6, 3 ]}{6}{9+9/17}{-62/17}

\ssi{17}{4}{[ 5, 2, 2, 2 ]}{13}{12+0}{-36/17}

\ssi{17}{5}{[ 4, 2, 3 ]}{7}{9+12/17}{-31/17}

\ssi{17}{10}{[ 2, 4, 2, 2 ]}{12}{11+5/17}{-24/17}

\ssi{18}{5}{[ 4, 3, 2 ]}{11}{9+8/9}{-2}

\ssi{18}{7}{[ 3, 3, 2, 2 ]}{13}{11+1/9}{-11/9}

\ssi{19}{3}{[ 7, 2, 2 ]}{13}{11+16/19}{-75/19}

\ssi{19}{4}{[ 5, 4 ]}{5}{9+9/19}{-68/19}

\ssi{19}{7}{[ 3, 4, 2 ]}{11}{9+18/19}{-39/19}

\ssi{19}{8}{[ 3, 2, 3, 2 ]}{12}{11+1/19}{-22/19}

\hline
\end{tabular}
$$

$$
\begin{tabular}{|c|c|c|c|c|}

\ssi{n}{q}{n/q = [b_1, \dots , b_s]}{q^{\prime}}{B_{\si{n}{q}}}{h_{\si{n}{q}}}

 & & & &   \\

\ssi{20}{3}{[ 7, 3 ]}{7}{10+1/2}{-23/5}

\ssi{20}{11}{[ 2, 6, 2 ]}{11}{11+1/10}{-16/5}

\ssi{21}{8}{[ 3, 3, 3 ]}{8}{9+16/21}{-13/7}

\ssi{21}{13}{[ 2, 3, 3, 2 ]}{13}{11+5/21}{-4/3}

\ssi{22}{5}{[ 5, 2, 3 ]}{9}{10+7/11}{-30/11}

\ssi{23}{3}{[ 8, 3 ]}{8}{11+11/23}{-128/23}

\ssi{23}{4}{[ 6, 4 ]}{6}{10+10/23}{-104/23}

\ssi{23}{5}{[ 5, 3, 2 ]}{14}{10+19/23}{-67/23}

\ssi{23}{7}{[ 4, 2, 2, 3 ]}{10}{11+17/23}{-42/23}

\ssi{24}{5}{[ 5, 5 ]}{5}{10+5/12}{-9/2}

\ssi{24}{7}{[ 4, 2, 4 ]}{7}{10+7/12}{-8/3}

\ssi{25}{7}{[ 4, 3, 2, 2 ]}{18}{12+0}{-52/25}

\ssi{25}{9}{[ 3, 5, 2 ]}{14}{10+23/25}{-3}

\ssi{26}{7}{[ 4, 4, 2 ]}{15}{10+11/13}{-38/13}

\ssi{27}{4}{[ 7, 4 ]}{7}{11+11/27}{-148/27}

\ssi{27}{5}{[ 6, 2, 3 ]}{11}{11+16/27}{-11/3}

\ssi{27}{8}{[ 4, 2, 3, 2 ]}{17}{11+25/27}{-2}

\ssi{28}{5}{[ 6, 3, 2 ]}{17}{11+11/14}{-27/7}

\ssi{29}{5}{[ 6, 5 ]}{6}{11+11/29}{-158/29}

\ssi{29}{8}{[ 4, 3, 3 ]}{11}{10+19/29}{-79/29}

\ssi{29}{12}{[ 3, 2, 4, 2 ]}{17}{12+0}{-60/29}

\ssi{30}{11}{[ 3, 4, 3 ]}{11}{10+11/15}{-14/5}

\ssi{31}{7}{[ 5, 2, 4 ]}{9}{11+16/31}{-111/31}

\ssi{31}{11}{[ 3, 6, 2 ]}{17}{11+28/31}{-123/31}

\ssi{31}{12}{[ 3, 3, 2, 3 ]}{13}{11+25/31}{-58/31}

\ssi{33}{7}{[ 5, 4, 2 ]}{19}{11+26/33}{-127/33}

\ssi{34}{9}{[ 4, 5, 2 ]}{19}{11+14/17}{-66/17}

\ssi{34}{13}{[ 3, 3, 3, 2 ]}{21}{12+0}{-35/17}

\ssi{37}{8}{[ 5, 3, 3 ]}{14}{11+22/37}{-135/37}

\ssi{39}{14}{[ 3, 5, 3 ]}{14}{11+28/39}{-49/13}

\ssi{40}{11}{[ 4, 3, 4 ]}{11}{11+11/20}{-18/5}

\ssi{41}{11}{[ 4, 4, 3 ]}{15}{11+26/41}{-151/41}

\hline
\end{tabular}
$$

\newpage

\section*{Appendix $B$} \label{appendix A}
This appendix contains the classification of finite groups of automorphisms
 acting on Riemann surfaces of genus $2$, and  $3$ so that the quotient is
 isomorphic to $\mathbb{P}^1$.
In the last case we listed only the nonabelian groups.
Tables \ref{2-abeliani}, \ref{2-nonabelian} and
 \ref{3-nonabelian} are adapted from [Br90, pages 252, 254, 255].
For every $G$ we give a presentation, the vector $\mathbf{m}$ of branching data  and
the \verb|IdSmallGroup|$(G)$, that is the number of $G$ in the \verb|GAP4|
database of small groups. The second author  wishes to thank
S. A. Broughton who kindly communicated to him that the group $G(48,33)$
(Table \ref{3-nonabelian}, case $(3t)$)
was missing in \cite{Br90}.

\begin{table}[ht!]
\begin{center}
\begin{tabular}{|c|c|c|c|}
\hline
$ $ & $ $ & \verb|IdSmall| & $ $ \\
$\textrm{Case}$ & $G$ & \verb|Group|$(G)$ &  $\mathbf{m}$ \\
\hline \hline
$(1a)$ & $\mZ_2$ & $G(2,1)$ & $(2^6)$ \\
\hline
$(1b)$ & $\mZ_3$ & $G(3,1)$ & $(3^4)$ \\
\hline
$(1c)$ & $\mZ_4$ & $G(4,1)$ & $(2^2, 4^2)$ \\
\hline
$(1d)$ & $\mZ_2 \times \mZ_2$ & $G(4,2)$ & $(2^5)$ \\
\hline
$(1e)$  & $\mZ_5$  & $G(5,1)$ & $(5^3)$\\
\hline
 $(1f)$  & $\mZ_{6}$ & $G(6,2)$ & $(2^2,3^2)$ \\
\hline
$(1g)$  & $\mZ_6$ & $G(6,2)$ & $(3,6^2)$ \\
\hline
$(1h)$  & $\mZ_8$ & $G(8,1)$ & $(2,8^2)$ \\
\hline
$(1i)$  & $\mZ_{10}$ & $G(10,2)$ & $(2,5,10)$ \\
\hline
$(1j)$  & $\mZ_2 \times \mZ_6$ & $G(12,5)$ & $(2,6^2)$ \\
\hline
\end{tabular}
\end{center}
\caption{Abelian automorphism groups with rational quotient on  Riemann
  surfaces of genus $2$}
\label{2-abeliani}
\end{table}

\begin{table}[ht!]
\begin{center}
\begin{tabular}{|c|c|c|c|c|}
\hline
$ $ & $ $ &  \verb|IdSmall| & $ $ & $ $ \\
${\textrm{Case}}$ & $G$ & \verb|Group|$(G)$ & $\mathbf{m}$ & $ \textrm{Presentation}$ \\
\hline \hline
 $(2a)$  & $\mathcal{S}_3$  & $G(6,1)$ & $(2^2,3^2)$ & $\langle x,y\;|\;x=(123), \; y=(12)
 \rangle$  \\
\hline
$$ & $$ & $$ & $$  & $\langle i,j,k\; | \; i^2=j^2=k^2=-1$,
\\ $(2b)$ & $Q_8$ & $G(8,4)$ & $(4^3)$ & $ij=k, \; jk=i, \; ki=j \rangle$ \\
\hline
$(2c)$ & $D_4$ & $G(8,3)$ & $(2^3,4)$ & $\langle  x,y\; | \; x^2=y^4=1,
 \; xyx^{-1}=y^{-1} \rangle$ \\
\hline
$(2d)$ & $D_{4,3,-1}$ & $G(12,1)$ & $(3,4^2)$ & $\langle x,y\; | \; x^4=y^3=1,
 \; xyx^{-1}=y^{-1} \rangle$ \\
\hline
$(2e)$ & $D_6$ & $G(12,4)$ & $(2^3,3)$ & $\langle  x,y \; | x^2=y^6=1,
 \; xyx^{-1}=y^{-1} \rangle$ \\
\hline
$(2f)$ & $D_{2,8,3}$ & $G(16,8)$ & $(2,4,8)$ & $\langle  x,y \; | \;x^2=y^8=1,
 \; xyx^{-1}=y^3 \rangle$ \\
\hline
$$ & $$ & $$ & $$ & $\langle x,y,z,w \; | \; x^2=y^2=z^2=w^3=1,$
\\
$(2g)$ & $G=\mathbb{Z}_2 \ltimes ((\mathbb{Z}_2)^2 \times
\mathbb{Z}_3)$ & $G(24,8)$ & $(2,4,6)$ & $[y,z]=[y,w]=[z,w]=1,$
\\
$$ & $$ & $$ & $$ & $xyx^{-1}=y,\; xzx^{-1}=zy, \; xwx^{-1}=w^{-1} \rangle$ \\
\hline
$(2h)$ &  $\textrm{SL}_2(\mathbb{F}_3)$ & $G(24,3)$ & $(3^2,4)$ & $ \langle x,y \; | \;
 x=
\left(
              \begin{array}{cc}
                1 & 1 \\
                0 & 1 \\
              \end{array}
            \right), \;
y=
\left(
              \begin{array}{cc}
                \;\;0 & \;\;1 \\
                -1 & -1 \\
              \end{array}
            \right)
\rangle$ \\
\hline
$(2i)$ &  $\textrm{GL}_2(\mathbb{F}_3)$ & $G(48,29)$ & $(2,3,8)$ & $ \langle x,y \; | \;
 x=
\left(
              \begin{array}{cc}
                1 & \;\;1 \\
                0 & -1 \\
              \end{array}
            \right), \;
 y=
\left(
              \begin{array}{cc}
                0 & -1 \\
                1 & -1 \\
              \end{array}
            \right)
\rangle$ \\
\hline
\end{tabular}
\end{center}
\caption{Nonabelian automorphism groups with rational quotient
  on Riemann surfaces of genus $2$.}
\label{2-nonabelian}
\end{table}

\newpage

\begin{table}[ht!]
\begin{center}
\begin{tabular}{|c|c|c|c|c|}
\hline
$ $ & $ $ &  \verb|IdSmall| & $ $ & $ $ \\
${\textrm{Case}}$ & $G$ & \verb|Group|$(G)$ & $\mathbf{m}$ & $ \textrm{Presentation}$ \\
\hline \hline $(3a)$  & $\mathcal{S}_3$  & $G(6,1)$ & $(2^4,3)$ & $\langle x,y\;|\;x=(12), \;
y=(123)
 \rangle$  \\
 \hline
$(3b)$ & $D_4$ & $G(8,3)$ & $(2^2,4^2)$ & $\langle  x,y\; | \; x^2=y^4=1,
 \; xyx^{-1}=y^{-1} \rangle$ \\
\hline $(3c)$ & $D_4$ & $G(8,3)$ & $(2^5)$ & $\langle  x,y\; | \; x^2=y^4=1,
 \; xyx^{-1}=y^{-1} \rangle $ \\
\hline
 $(3d)$ & $D_{4,3,-1}$ & $G(12,1)$ & $(4^2,6)$ & $\langle x,y\; | \; x^4=y^3=1,
 \; xyx^{-1}=y^{-1} \rangle$ \\
\hline $(3e)$ & $D_6$ & $G(12,4)$ & $(2^3,6)$ & $\langle  x,y\; | \; x^2=y^6=1,
 \; xyx^{-1}=y^{-1} \rangle$ \\
 \hline
$(3f)$ & $\mathcal{A}_4$ & $G(12,3)$ & $(2^2, 3^2)$ & $\langle x,y \; | \;
x=(12)(34), \; y=(123) \rangle$ \\
\hline $(3g)$ & $D_{2,8,5}$ & $G(16,6)$ & $(2,8^2)$ & $\langle x,y\; | \; x^2=y^8=1,
 \; xyx^{-1}=y^5 \rangle$ \\
\hline $(3h)$ & $D_{4,4,-1}$ & $G(16,4)$ & $(4^3)$ & $\langle x,y\; | \; x^4=y^4=1,
 \; xyx^{-1}=y^{-1} \rangle$ \\
\hline $(3i)$ & $\mathbb{Z}_2 \times D_4$ & $G(16,11)$ & $(2^3,4)$ & $\langle z \; | \;
z^2=1 \rangle \times   \langle  x,y\; | \; x^2=y^4=1,
 \; xyx^{-1}=y^{-1} \rangle $ \\
\hline $ $ & $ $ & $ $ & $ $ & $\langle x,y,z \; | \; x^2=y^2=z^4=1,$\\
$(3j)$ & $\mathbb{Z}_2 \ltimes (\mathbb{Z}_2 \times
\mathbb{Z}_4) $ & $G(16,13)$ & $(2^3,4)$ & $[x,z]=[y,z]=1, \; xyx^{-1}=yz^2 \rangle$ \\
\hline $(3k)$ & $D_{3,7,2}$ & $G(21,1)$ & $(3^2,7)$ & $\langle x,y\; | \; x^3=y^7=1,
 \; xyx^{-1}=y^2 \rangle$ \\
\hline $(3l)$ & $D_{2,12,5}$ & $G(24,5)$ & $(2,4,12)$ & $\langle x,y\; | \; x^2=y^{12}=1,
 \; xyx^{-1}=y^5 \rangle$ \\
\hline $(3m)$ & $\mathbb{Z}_2 \times \mathcal{A}_4$ & $G(24,13)$ & $(2,6^2)$ & $ \langle z \; | \;
z^2=1 \rangle \times \langle x,y \; | \; x=(12)(34), \; y=(123) \rangle
$ \\
\hline $(3n)$ & $\textrm{SL}_2(\mathbb{F}_3)$ & $G(24,3)$ & $(3^2,6)$ & $ \langle x,y \;
| \;
 x=
\left(
              \begin{array}{cc}
                1 & 1 \\
                0 & 1 \\
              \end{array}
            \right), \;
y= \left(
              \begin{array}{cc}
                \;\;0 & \;\;1 \\
                -1 & -1 \\
              \end{array}
            \right)
\rangle$ \\
\hline
$(3o)$ & $\mathcal{S}_4$ & $G(24,12)$ & $(3,4^2)$ & $\langle x,y \; | \; x=(1234), \; y=(12)   \rangle$ \\
\hline
$(3p)$ & $\mathcal{S}_4$ & $G(24,12)$ & $(2^3,3)$ & $\langle x,y \; | \; x=(1234), \; y=(12)   \rangle$ \\
\hline $ $ & $ $ & $ $ & $ $ & $\langle x,y,z \; | \; x^2=y^2=z^8=1,$\\
$(3q)$ & $\mathbb{Z}_2 \ltimes(\mathbb{Z}_2 \times
\mathbb{Z}_8)$ & $G(32,9)$ & $(2,4,8)$ & $[x,y]=[y,z]=1, \; xzx^{-1}=yz^3 \rangle$ \\
\hline
$ $ & $ $ & $ $ & $ $ & $\langle x,y,z \; | \; x^2=y^2=z^8=1,$\\
$(3r)$ & $\mathbb{Z}_2 \ltimes D_{2,8,5}$ & $G(32,11)$ & $(2,4,8)$ & $yzy^{-1}=z^5, \; xyx^{-1}=yz^4, \; xzx^{-1}=yz^3 \rangle$ \\
\hline $(3s)$ & $\mathbb{Z}_2 \times \mathcal{S}_4$ & $G(48,48)$ & $(2,4,6)$ & $ \langle z \; | \;
z^2=1 \rangle \times \langle x,y \; | \; x=(12), \; y=(1234) \rangle
$ \\
\hline $ $ & $ $ & $ $ & $ $ & $ \langle x,y,z,w,t \; | \; x^2=z^2=w^2=t,
\; y^3=1,\; t^2=1, $ \\
$(3t)$ & $G(48,33)$ & $G(48,33)$ & $(2,3,12)$ & $yzy^{-1}=w, \; ywy^{-1}=zw, \; zwz^{-1}=wt$,\\
$ $ & $ $ & $ $ & $ $ & $[x,y]=[x,z]=1 \rangle$ \\   \hline
$ $ & $ $ & $ $ & $ $ & $\langle x,y,z \; | \; x^3=y^4=z^4=1,$\\
$(3u)$ & $\mathbb{Z}_3 \ltimes(\mathbb{Z}_4)^2$ & $G(48,3)$ &
$(3^2,4)$ & $[y,z]=1, \; xyx^{-1}=z, \; xzx^{-1}=(yz)^{-1} \rangle$ \\
\hline
$$ & $$ & $$ & $$ & $\langle x,y,z,w \; | \; x^2=y^3=z^4=w^4=1,$
\\
$(3v)$ & $\mathcal{S}_3 \ltimes (\mathbb{Z}_4)^2$ & $G(96,64)$ & $(2,3,8)$ & $[z,w]=1, \;
xyx^{-1}=y^{-1}, \; xzx^{-1}=w,$
\\
$$ & $$ & $$ & $$ & $xwx^{-1}=z,\; yzy^{-1}=w, \; ywy^{-1}=(zw)^{-1} \rangle$ \\\hline
$(3w)$ & $\textrm{PSL}_2(\mathbb{F}_7)$ & $G(168,42)$ & $(2,3,7)$ &
$\langle x,y \; | \; x=(375)(486), \; y=(126)(348)\rangle$ \\
\hline
\end{tabular}
\end{center}
\caption{Nonabelian automorphism groups with rational quotient
  on Riemann surfaces of genus $3$.}
\label{3-nonabelian}
\end{table}

\bigskip \bigskip
ERNESTO MISTRETTA \\
Università degli Studi di Padova,
Dipartimento di Matematica Pura e Applicata,
Via Trieste 63,
35121 Padova, ITALY  \\
\emph{E-mail address}: \verb|ernesto@math.unipd.it| \\ \\

FRANCESCO POLIZZI \\
Dipartimento di Matematica,
Università della Calabria,
Via P. Bucci Cubo 30B,
87036 Arcavacata di Rende (CS), Italy. \\
\emph{E-mail address}: \verb|polizzi@mat.unical.it|


\newpage

\begin{thebibliography}{999999}










\bibitem[BPV84]{BPV}
W. Barth, C. Peters, A. Van de Ven: \emph{Compact Complex Surfaces}, Springer-Verlag
1984.

\bibitem[Bar99]{Bar99}
R. Barlow: Zero-cycles on Mumford's surface, \emph{Math. Proc. Camb. Phil. Soc.}
$\boldsymbol{126}$ (1999), 505-510.




\bibitem[Be96]{Be2}
A. Beauville: \emph{Complex algebraic surfaces}, Cambridge University Press 1996.




\bibitem[Bre00]{Bre00}
T. Breuer: \emph{Characters and Automorphism groups of Compact Riemann Surfaces},
 Cambridge University Press 2000.

\bibitem[Br90]{Br90}
S. A. Broughton: Classifying finite group actions on surfaces of low genus, \emph{J. of
Pure and Applied Algebra} $\boldsymbol{69}$ (1990), 233-270.

\bibitem[Ca81]{Ca}
F. Catanese: On a class of surfaces of general type, in \emph{Algebraic Surfaces}, CIME,
Liguori (1981), 269-284.



\bibitem[Ca99]{Ca99}
F. Catanese: Singular bidouble covers and the construction of interesting algebraic
surfaces, \emph{Contemporary Mathematics} $\boldsymbol{241}$ (1999), 97-119.





\bibitem[CaCi91]{CaCi91}
F. Catanese and C. Ciliberto: Surfaces with $p_g=q=1$, \emph{Symposia Math.}
$\boldsymbol{32}$ (1991), 49-79.

\bibitem[CaCi93]{CaCi93}
F. Catanese and C. Ciliberto: Symmetric product of elliptic curves and surfaces of
general type with $p_g=q=1$, \emph{Journal of Algebraic Geometry} $\boldsymbol{2}$
(1993), 389-411.


\bibitem[CaPi06]{CaPi06}
F. Catanese, R. Pignatelli: Fibrations of low genus I, \emph{Ann.
Sci. École Norm. Sup.} (4) $\boldsymbol{39}$ (2006), no. 6,
1011-1049.

\bibitem[CarPol]{CarPol}
G. Carnovale, F. Polizzi: The classification of surfaces of general
type with $p_g=q=1$ isogenous to a product, e-print
$\mathbf{arXiv:0704.0446}$ (2007), to appear in \emph{Adv. Geom.}












\bibitem[FK92]{FK92}
H. M. Farkas, I. Kra: \emph{Riemann Surfaces}, Graduate Texts in Mathematics
$\boldsymbol{71}$, $2^{nd}$ Edition, Springer-Verlag 1992.



\bibitem[Fre71]{Fre71}
E. Freitag: Uber die Struktur der Funktionenk$\ddot{\textrm{o}}$rper zu hyperabelschen
Gruppen I, \emph{J. Reine. Angew. Math.} $\boldsymbol{247}$ (1971), 97-117.




\bibitem[GAP4]{GAP4}
The GAP~Group, \emph{GAP -- Groups, Algorithms, and Programming,
  Version 4.4}; 2006, http://www.gap-system.org.







\bibitem[H71]{H71} W. J. Harvey: On the branch loci in
Teichm{\"u}ller space, \emph{Trans. Amer. Mat. Soc.} $\boldsymbol{153}$ (1971), 387-399.








\bibitem[Is05]{Is05}
H. Ishida: Bounds for the relative Euler-Poincaré characteristic of certain hyperelliptic
fibrations, \emph{Manuscripta Math.} $\boldsymbol{118}$ (2005), 467-483.

\bibitem[Is06]{Is06}
H. Ishida: Catanese-Ciliberto surfaces of fiber genus three with unique singular fibre,
\emph{Tohoku Math. J.} $\boldsymbol{58}$ (2006), 33-69.











\bibitem[Lau71]{Lau71}
H. B. Laufer: \emph{Normal two-dimensional singularities}, Annals of Mathematics Studies $\boldsymbol{71}$,
 Princeton University Press 1971.













\bibitem[Pi08]{Pi08}
R. Pignatelli: Some (big) irreducible components of the moduli space
of minimal surfaces of general type with $p_g=q=1$ and $K^2=4$,
e-print $\mathbf{arXiv:0801.1112}$ (2008).



\bibitem[Pol09]{Pol09}
F. Polizzi: Standard isotrivial fibrations with $p_g=q=1$, \emph{J.
Algebra} $\boldsymbol{321}$ (2009), 1600-1631.


\bibitem[Ri07]{Ri07}
C. Rito: On surfaces with $p_g=q=1$ and non-ruled bicanonical
involution,  \emph{Ann. Sc. Norm. Super. Pisa Cl. Sci} (5)  6
(2007), no. 1, 81-102.

\bibitem[Ri08]{Ri08}
C. Rito: On equations of double planes with $p_g=q=1$, e-print
$\mathbf{arXiv:0804.2227}$ (2008).








\bibitem[Se96]{Se96}
F. Serrano: Isotrivial fibred surfaces, \emph{Annali di Matematica pura e applicata},
vol. CLXXI (1996), 63-81.









\end{thebibliography}
\end{document}